\documentclass[preprint]{elsarticle}
\usepackage{amssymb}
\usepackage{amsmath}
\newcommand{\R}{\mathbb{R}}

\newcommand{\C}{\mathbb{C}}
\newcommand{\Z}{\mathbb{Z}}
\font\eufm=eufm10
\def\frak#1{\hbox{\eufm#1}}
\newcommand{\bd}{\begin{document}}
\newcommand{\ed}{\end{document}}
\newcommand{\be}{\begin{enumerate}}
\newcommand{\ee}{\end{enumerate}}
\newcommand{\bi}{\begin{itemize}}
\newcommand{\ei}{\end{itemize}}
\newcommand{\ba}{\begin{array}}
\newcommand{\ea}{\end{array}}
\newcommand{\vs}{\vspace*{0.3\baselineskip}}%%%maly odstep pionowy
\newcommand{\vsm}{\vspace*{-0.3\baselineskip}}%%%maly odstep pionowy ujemny
\newcommand{\kom}[1]{{\em #1}\newline}%%%komentarz w oddzielnej lini
\newtheorem{defi}{Definition}[section]
\newtheorem{tw}[defi]{Theorem}
\newtheorem{prop}[defi]{Proposition}
\newtheorem{lem}[defi]{Lemma}
\newtheorem{re}[defi]{Remark}
\newtheorem{col}[defi]{Corollary}
\newtheorem{ex}[defi]{Examples}
\newtheorem{zad}{Exercise}[section]
\newtheorem{zal}{Assumptions}[section]
\newtheorem{assumpt}[defi]{Assumptions}
\newcommand{\Om}{\Omega}
\newcommand{\om}{\omega}
\newcommand{\G}{\Gamma}
\newcommand{\D}{\Delta}
\renewcommand{\d}{\delta}
\newcommand{\ga}{\gamma}
\newcommand{\eps}{\epsilon}
\newcommand{\ove}{\overline}
\newcommand{\ms}{\oplus}
\newcommand{\mt}{\otimes}
\newcommand{\dz}{\wedge}
\newcommand{\lra}{\longrightarrow}
\newcommand{\rra}{\rightrightarrows}
\newcommand{\rel}{\mbox{$\,$\rule[0.5ex]{1.1em}{0.2pt}$\triangleright\,$}}
\newcommand{\dow}{\hspace*{\fill}\rule{1.6ex}{1.6ex}\hspace*{1em}}
\newcommand{\dowl}{\hspace*{\fill}\rule{1ex}{1ex}\hspace*{1em}}
\newcommand{\sd}{\hspace{0.3ex}\tiny{\rhd\mbox{\hspace{-2ex}}<}\hspace{0.3ex}}
\newcommand{\mmt}[2]{\mbox{$\vphantom{}_{#1}\times_{#2}$}}
\newcommand{\g}{\frak g}
\newcommand{\ab}{\frak a}
\newcommand{\bb}{\frak b}
\newcommand{\h}{\frak h}
\newcommand{\got}{\frak t}
\newcommand{\gotg}{\frak g}
\newcommand{\gota}{\frak a}
\newcommand{\gotb}{\frak b}
\newcommand{\gotc}{\frak c}
\newcommand{\gothh}{\frak h}
\newcommand{\gott}{\frak t}
\newcommand{\hd}{\hat{\d}}
\newcommand{\oml}{\Omega_L^{1/2}}
\newcommand{\omr}{\Omega_R^{1/2}}
\newcommand{\omh}{\Omega^{1/2}}
\newcommand{\lo}{\lambda_0}
\newcommand{\ro}{\rho_0}
\newcommand{\sA}{\mbox{$\mathcal A$}}
\newcommand{\sT}{\mbox{$\mathcal T$}}
\newcommand{\sB}{\mbox{$\mathcal B$}}
\newcommand{\sF}{\mbox{$\mathcal F$}}
\newcommand{\sO}{\mbox{$\mathcal O$}}
\newcommand{\sD}{\mbox{$\mathcal D$}}
\newcommand{\lma}{\Lambda^{max}}
\newcommand{\timh}{\times_h}
\newcommand{\Gd}{\G^{(2)}}
\newcommand{\el}{e_L}
\newcommand{\er}{e_R}
\newcommand{\GG}{\G_1\times\G_2}
\newcommand{\gdot}{\hspace{-0.1em}\cdot\hspace{-0.1em}}
\newcommand{\tran}{\frown\hspace{-2.2ex}|\hspace{1.9ex}}
\newcommand{\la}[2]{\Lambda_{#1#2}}
\newcommand{\kad}{ad^{\#}}
\newcommand{\wl}[1]{\vphantom{X}_{#1}{\G}}
\newcommand{\te}{\tilde{e}}
\newcommand{\notka}[1]{}
\newcommand{\dif}{differential }
\newcommand{\gru}{groupoid }
\newcommand{\grus}{groupoids }
\newcommand{\ti}{\tilde}
\newcommand{\halden}{half density }
\newcommand{\haldens}{half densities }
\renewcommand{\top}{topological }
\newcommand{\Setrel}{\mbox{\rm SetRel}}
\newcommand{\hX}{\mbox{$\hat{X}$}}
\newcommand{\cstardwa}{\mbox{$C^*_r(\Gamma\times\Gamma)$}}
\begin{document}
\title{On the quantum 'ax+b' group}
\author{Piotr Stachura}
\address{Faculty of Applied Informatics and Mathematics, Warsaw University of Life Sciences-SGGW,
ul Nowoursynowska 166, 02-787 Warszawa, Poland,  
e-mail: stachura@fuw.edu.pl}
\date{}
\begin{abstract} The more detailed description of the 'ax+b' group of Baaj and Skandalis is presented. 
In particular we give generators and present formulae for the action of the  comultiplication on them. 
We prove that this quantum group is defined by a  twist.
\end{abstract}
\maketitle
%%%%%%%%%%%%%%%%%%%%%%%%%%%%%%%%%%%%%%%%%%%%%%%%%%%%%%%%%%%%%%%%%%%%%%%%%%%%%%%%
%%%%%%%%%%%%%%%%%%%%%%%%%%%%%%%%%%%%%%%%%%%%%%%%%%%%%%%%%%%%%%%%%%%%%%%%%%%%%%%
\normalmarginpar
\section{Introduction}
The purpose of this work is to give a more detailed  description of the  quantum 'ax+b' group of Baaj and Skandalis 
\cite{BS, VV, BV}.  In particular we describe generators of the $C^*$-algebra, show that they satisfy relations which are 
formally equivalent to that given in \cite{VV}, prove that the comultiplication is given by a twist  and 
compute an action of the comultiplication on generators. 
We also show that this quantum group is a quantization of a Poisson-Lie structure on the  classsical 'ax+b' group. 

Strictly speaking {\em it is not proven} here that the quantum group constructed in this paper and the Baaj-Skandalis example coincide. 
They have the same underlying $C^*$-algebra and our comultiplication formally coincides with their on generators. 
Had the quantum ``ax+b'' been defined by a global decomposition, the results of \cite{DLG} would guarantee that these groups are the same. 
But this is not the case, so formally there is still (rather improbable)  possibility that the group presented here is not the Baaj-Skandalis example.

On the Hopf $*$-algebra level this group is given by generators $A\,,\,A^{-1}\,,\,Z$ and relations \cite{VV}:
$$A=A^*\,\,,\,\,\,Z=-Z^*\,,\,\, [A,Z]=A(1-A)$$ 
Together with  comultiplication $\Delta$, counit $\epsilon$ and antipode $S$:
$$\Delta(A)=A\mt A\,\,,\,\,\Delta(Z)= Z\mt A+ 1\mt Z\,\,,\,
\epsilon(A)=1\,,\,\,\epsilon(Z)=0\,\,,\,\,S(A)=A^{-1}\,\,,\,\,S(Z)=- Z A^{-1}$$
If one wants to find a $C^*$-algebra with affiliated elements $A,Z$ that satisfy these commutation relations, 
first step is to represent them on a Hilbert space. One hopes to find $A,A^{-1}$--unbounded, invertible, selfadjoint operators and
$iZ$--selfadjoint,  satisfying in some reasonable sense the relation  $[A,Z]=A(1-A)$. 
This is possible but the most natural choice leads to the situation where   $\Delta(iZ)$ is symmetric but not selfadjoint.

One can rewrite these relations using different generators and this is the way we choose.
Let us define $Y:=A^{-1}-I$ and $X:=i Z$. These elements satisfy:
\notka{relkom}

\begin{align}\label{relkom}
Y=Y^*\,,\,\,X=X^*\,,\,\, [X,Y]=i Y\,,\nonumber\\
\Delta(Y)=Y\mt Y+I\mt Y +Y\mt I\,,\,\,\Delta(X)= X\mt (Y+I)^{-1} + I\mt X
\end{align}
%%%%%%%
These are relations we are going to give a meaning to.

In the following we will use groupoid algebras, so now we recall basic facts and establish the 
relevant notation (see \cite{DLG,DG} for a detailed exposition). The category of groupoids used here is described in \cite{SZ1} and in a 
differential setting in \cite{SZ2}. {\em A groupoid} is a set $\Gamma$ together with a subset $E\subset\Gamma$
(the set of identities), an associative relation (multiplication) $m:\Gamma\times\Gamma\rel \Gamma$ 
and an involutive mapping (inverse)
$s:\Gamma\rightarrow\Gamma$. They satisfy certain relations that entail the existence of two projections 
$\el,\er:\Gamma\rightarrow E$ (target and source projections) and the fact that $m$ is a mapping defined on the set 
$\{(x,y)\in\Gamma\times\Gamma : \er(x)=\el(y)\}$ of composable pairs; this is the standard notion of groupoid. 
%{\em groupoid is a small category , such that every morphism is an isomorphism}.
 The definition of  {\em a morphism of groupoids} used here is different, however. 
By a morphism of groupoids $\Gamma, \Gamma'$  we mean  a relation $h:\Gamma\rel\Gamma'$ 
which  satisfies:  $hm=m'(h\times h)$, $h s=s' h$ and $hE=E'$ (see \cite{SZ1,SZ2}). 
In a differential setting $\G$ is a smooth(Hausdorff) 
manifold, the set of units and the set of composable pairs are closed submanifolds, $s$ is 
a diffeomorphism, $m$ is a differential reduction  and $\el,\er$ are surjective submersions. 
%For a composable pair 
%$\gamma,\gamma'$ we will usually write $\gamma \gamma'$ instead of $m(\gamma,\gamma')$. 
%For two subsets $A,B\subset\Gamma$ we write $AB:=\{\gamma_1 \gamma_2 : \gamma_1\in A, \gamma_2\in B\}$.

{\em A bisection} of a differential groupoid $\Gamma$ is a submanifold $B\subset \Gamma$ such that 
$\el|_B$ and $\er|_B$ are diffeomorphisms $B\rightarrow E$. If $B\subset \Gamma$ is a bisection and 
$h:\Gamma\rel\Gamma'$ is a morphism then $h(B)$ is a bisection of $\Gamma'$.

Let $\Gamma$ be a differential groupoid and  let $\Om_L^{1/2}, \Om_R^{1/2}$ denote the bundles of complex  half 
densities along left and right fibers. {\em A groupoid *-algebra} $\sA(\Gamma)$ is a vector space of compactly 
supported, smooth sections of $\Om_L^{1/2}\mt\Om_R^{1/2}$ together with a convolution and $*$-operation.
To write explicit  formulae let us choose $\lambda_0$ - a real, nonvanishing, left invariant half density along left fibers 
(in fact that means we choose a Haar system on $\Gamma$, the ``choice-free'' definition 
and other detailes are  given in \cite{DG}.) Let $\rho_0:=s(\lambda_0)$ be the corresponding right invariant half density, 
and $\om_0:=\lambda_0\mt \rho_0$. Then any $\om\in\sA(\Gamma)$ can be written as $\om=f\om_0$ for a function 
$f\in\sD(\Gamma)$ (smooth, compactly supported). With such a choice  we write $(f_1\om_0)\,(f_2\om_0)=:(f_1*f_2)\om_0$,
$(f\om_0)^*=:(f^*)\om_0$ and: \notka{mult}

\begin{equation}\label{mult}
(f_1*f_2)(\gamma):=\int \lambda_0^2(\gamma') f_1(\gamma')f_2(s(\gamma') \gamma)=
\int\rho_0^2(\gamma')f_1(\gamma s(\gamma')f_2(\gamma')\,\,,\,\,f^*(\gamma):=\overline{f(s(\gamma))}
\end{equation}
The first integral is over the left fiber passing through $\gamma$, the second over the right one.\\
The choice of $\om_0$ defines a norm that makes  $\sA(\Gamma)$ a normed $*$-algebra.:
$$||f\om_0||_0=:||f||_0=max\left\{\sup_{e\in E}\int_{\el^{-1}(e)}\lambda_0^2(\gamma)|f(\gamma)|,\,
\sup_{e\in E}\int_{\er^{-1}(e)}\rho_0^2(\gamma)|f(\gamma)|\right\}$$

There is a faithful representation $\pi_{id}$ of $\sA(\Gamma)$ on $L^2(\Gamma)$ described as follows: choose $\nu_0$ - a real, 
nonvanishing half density on $E$; since $\er$ is a surjective submersion one can define $\psi_0:=\rho_0\mt\nu_0$ - this is 
a real, nonvanishing, half density on $\Gamma$. For $\psi=f_2\psi_0,\, f_2\in\sD(\Gamma)$  the representation is given by 
$\pi_{id}(f_1\om_0) (f_2\psi_0)=:(\pi_{id}(f_1)f_2)\psi_0$ and $\pi_{id}(f_1)f_2=f_1*f_2$  is as in (\ref{mult}). 
The estimate $||\pi_{id}(\om)||\leq||\om||_0$ makes  possible the definition:
{\em  The reduced $C^*$-algebra of a groupoid} is the completion of $\sA(\Gamma)$ in the  norm $||\om||:=||\pi_{id}(\om)||$.
We will also use the following fact which is a direct consequence of the definition of the norm $||f||_0$:
\begin{lem}\label{ind-lim}\notka{ind-lim} Let $U\subset \Gamma$ be an open set with compact closure. There exists $M$ such that
$||f||_0\leq M \sup |f(\gamma)|$ for any $f\in\sD(\Gamma)$ with support in $U$. If $f_n\in\sD(\Gamma) $ 
have supports in a fixed compact set and $f_n$ converges to $f\in\sD(\Gamma)$ uniformly then $f_n\om_0$ converges 
to $f\om_0$ in  $C^*_r(\Gamma)$.
\dowl
\end{lem}

A morphism $h:\Gamma\rel\Gamma'$ of differential groupoids defines a mapping 
$\hat{h}:\sA(\Gamma)\rightarrow L(\sA(\Gamma'))$ (linear mappings of $\sA(\Gamma')$), 
this mapping commutes with (right) multiplication in $\sA(\Gamma')$ i.e. 
$$\hat{h}(\om)(\om')\om''=\hat{h}(\om)(\om' \om'')\,\,,\,\om\in\sA(\Gamma)\,,\,\om',\om''\in\sA(\Gamma')$$
and we use  notation $\hat{h}(\om)\om'$ (see formula (\ref{delta0-form}) in the Appendix for an example of such a mapping);
 there is also a representation $\pi_h$ of $\sA(\Gamma)$ on 
$L^2(\Gamma')$; these objects satisfy some obvious compatibility conditions with respect to multiplication 
and $*$-operation (see \cite{DG} for a detailed exposition).

Now we recall some facts about {\em double groups}. Let $G$ be a group and $A,B\subset G$ subgroups 
such that $A\cap B=\{e\}$.  Every element $g$ in the set   $\Gamma:=AB\cap BA$ can be written uniquely as 
$$g=a_L(g) b_R(g)=b_L(g) a_R(g)\,;\,a_L(g),  a_R(g)\in A\,,\,b_L(g),b_R(g)\in B.$$
These decompostions define surjections: $a_L,a_R: \Gamma\rightarrow A$ and 
$b_L,b_R: \Gamma\rightarrow B$  (in fact $a_L, b_R$ are defined on $AB$ and $b_L, a_R$ on $BA$, 
we will denote these extensions by the same symbols). The formulae:
\notka{dlg-p}\begin{equation}\label{dlg-p}E:=A\,,\,\,s(g):=b_L(g)^{-1}a_L(g)=a_R(g) b_R(g)^{-1}\,,\,
Gr(m_A):=\{(b_1 a b_2; b_1 a, a b_2) :b_1a, ab_2\in\Gamma\}\end{equation}
define the structure of the  groupoid $\Gamma_A$ over $A$ on $\Gamma$. 
The analogous formulae define the groupoid  $\Gamma_B$ over $B$. 
On the other hand for a subgroup $B\subset G$ there is a (right) transformation groupoid 
$(B\setminus G)\times B$. The following lemma explains relation between these groupoids.
\notka{embed}\begin{lem}\label{embed}\notka{embed} The map:
$$\Phi:\Gamma_A\ni g\mapsto ([a_L(g)], b_R(g))\in (B\backslash G)\times B$$
is an isomorphism of the groupoid $\Gamma_A$ over $A$ with the restriction of a 
(right) transformation groupoid $ (B\backslash G)\times B$ to the set $\{[a] : a\in A\}\subset B\backslash G$.
\end{lem}
{\em Proof:} We give only the sketch of the proof. Let $\G_A\ni g=ab=b'a'$. Then $\Phi(g)=([a],b)$ and $[a]\cdot b=[ab]=[a']$, 
so $\Phi(\G_A)$ is really contained in  the restriction. 
On the other hand, if $([g],b)$ is an element of  the restriction i.e $[g]=[a]$ and $[gb]=[a']$ then the
 mapping $([g],b)\mapsto a_R(g)b\in G$ is well defined, has image in $\G_A$ and is the inverse of $\Phi$. Direct computations show that $\Phi$ 
is an isomorphism of groupoids. 
\dowl

If $AB=G$ (i.e. $\Gamma=G$)  the  triple $(G;A,B)$ is called {\em a double group} 
and in this situation we will denote groupoids 
$\Gamma_A, \Gamma_B$ by $G_A, G_B$. It turns out that $\delta:=m_B^T: G_A\rel G_A\times G_A$ is a coassociative morphism of groupoids. 
Applying the lemma \ref{embed} to the groupoid $G_A$ we can identify it with the transformation groupoid 
$(B\backslash G)\times B$. So $G_A=A\times B$ is a right transformation groupoid for the action
$(a,b)\mapsto a_R(ab)$ i.e the structure is given by:
$$E:=\{(a,e): a\in A\}\,,\,s(a,b):=(a_R(ab),b^{-1}),$$
$$m:=\{(a_1,b_1 b_2;a_1, b_1,a_R(a_1 b_1),b_2): a_1\in A\,,\,b_1,b_2\in B\}$$
{\em In the formula above, we identify a relation $m:\Gamma\times \Gamma\rel\Gamma$ with its graph, i.e. subset of 
$\Gamma\times \Gamma\times \Gamma$. We will use such notation throughout the paper.}
If $G$ is a Lie group, $A,B$ are closed sugroups, $A\cap B=\{e\}\,,\,AB=G$ 
then $(G;A,B)$ is called {\em a double Lie group.} 
It turns out that the mapping $\hat{\delta}$, defined by the morphism $\delta$ (compare (\ref{delta0-form}) in Appendix),
extends to the coassociative morphism $\Delta$ of $C^*_r(G_A)$ and $C^*_r(G_A\times G_A)=C^*_r(G_A)\mt C^*_r(G_A)$ 
which satisfies {\em density conditions} ({\em cls} denotes the closed linear span):
$$cls\{\Delta(a)(I\mt b) :a,b \in C^*_r(G_A)\}=cls\{\Delta(a)(b\mt I) :a,b \in C^*_r(G_A)\}=C^*_r(G_A)\mt C^*_r(G_A)$$
There are other objects that make the pair $(C^*_r(G_A), \Delta)$ a locally compact quantum group, we refer to  \cite{DLG} for details.

The quantum  'ax+b' does not completely fit  into this framework but, as we will see, it is possible to describe it 
using this approach as a guiding line. The main technical problem is that vector fields that ``should'' define 
operators affiliated to a groupoid $C^*$-algebra are not complete, so corresponding  operators are not essentialy 
selfadjoint on their ``natural'' domains. So one has to choose right domains or overcome this problem in a different way.

The next section ``sets the stage''; in the third one we consider a general situation in which a twist can be defined and apply results
 to the 'ax+b' group in the fourth section. (The situation from the third section also appears in the $\kappa$-Poincare Group, 
this will be decribed in a forthcoming paper). In the fifth section we give generators, show that they satisfy relations 
formally equivalent to ones given in the beginning of this introduction  and compute the action of 
comultiplication on them; we also express the twist as a function of generators. In the last section we consider 
our group as a deformation of a Poisson-Lie group. In the appendix we collect some formulae used in the paper.

As pointed out by the Referee, there appeared a preprint \cite{tuset} where the Baaj-Skandalis example is supposed 
to be given by a twist (this fact is not proven there). That work seems to be more general but uses completely different approach.

%%%%%%%%%%%%%%%%%%%%%%%%%%%%%%%%%%%%%%%%%%%%%%%%%%%%%%%%%%%%%%%%%%%%%%%%%%%%%%%%%%%%%%%%%%%%%%%%%%%%%%%%%%%%%%%%%%%
%%%%%%%%%%%%%%%%%%%%%%%%%%%%%%%%%%%%%%%%%%%%%%%%%%%%%%%%%%%%%%%%%%%%%%%%%%%%%%%%%%%%%%%%%%%%%%%%%%%%%%%%%%%%%%%%%%
\section{Setup}\label{setup}

Let $G$ be the 'ax+b' group i.e.
$G:=\R\times \R_*$ with the multiplication 
$(b_1,a_1)(b_2,a_2):=(b_1+a_1 b_2,a_1 a_2)$. 
%Let us consider the following subgroups of $G$:
For $s\in \R\cup\{\infty\}$ let us define the closed subgroup of $G$:
$$C_s:=\{(b,1+s b): 1+sb\neq 0\}\,,\,s\in\R\,;\,C_\infty:=\{(0,a) :  a\in \R_*\}.$$
(If we treat $\R\cup\{\infty\}$ as a one point compactification of $\R$ the mapping $s\mapsto C_s$ is continuous 
as a mapping into closed subgroups of $G$ with Fell topology).
Let $B:=C_0$ and $A:=C_\infty$.
Note that for $s\neq t$ we have  $C_s\cap C_t=\{e\}$ and for $s\in\R_*$:
$C_s=\{(\frac{c-1}{s}, c): c\in\R_*\}$.

For $s\neq t$ let:
$$\Gamma_{st}:=C_s C_t\cap C_t C_s=\{(b,a) : s t b-s a +t\neq 0\,,\,s t b- t a +s\neq 0\}=\Gamma_{ts}$$
It is straightforward to check that each $\Gamma_{st}$ is open and dense in $G$, and for $t\neq 0$ $\Gamma_{0t}=G$;
$\Gamma_{st}$ is a differential groupoid over $C_s$ and $C_t$. 

$G$ acts on $\{C_s: s\in\R\cup\{\infty\}\}$ by adjoint action and 
this action induces isomorphisms of corresponding differential groupoids. So it is sufficient to consider the family 
$$\Gamma_s:=\Gamma_{s\infty}=C_s A\cap A C_s=\{(b,a): (1+s b)(a-s b)\neq 0\}$$
and grupoid structures over $C_s$ and $A$. Projections in $\Gamma_s$ on $C_s$ and $A$ will be denoted by 
$\tilde{a}_L,\tilde{a}_R:\Gamma_s\rightarrow A$ 
and $\tilde{c}_L,\tilde{c}_R: \Gamma_s\rightarrow C_s$, they are given by:
\begin{eqnarray} \tilde{a}_L(b,a):=(0,a-sb)\,,\,\,\tilde{a}_R(b,a):=(0,\frac{a}{1+sb}),\nonumber\\
\tilde{c}_L(b,a):=(b,1+sb)\,,\,\,\tilde{c}_R(b,a):=(\frac{b}{a-s b},\frac{a}{a-s b})\nonumber
\end{eqnarray}
Remaining parts  of structures of groupoids are as follows.\\
For $\Gamma_s\rra C_s$: 
the inverse $\tilde{s}_C(b,a):=(\frac{b}{a- s b},\frac{1+ s b}{a- s b})$ 
and the multiplication relation:
\begin{equation}\label{mc}\tilde{m}_C:=\{(b_1, \frac{a_1 a_2}{1+s b_2};b_1,a_1,b_2,a_2)\,:\,b_1=b_2 (a_1-s b_1)\}
\subset\Gamma_s\times\Gamma_s\times\Gamma_s\end{equation}
For $\Gamma_s\rra A$: 
the inverse $\tilde{s}_A(b,a):=(\frac{-b}{1+s b},\frac{a- s b}{1+s b})$  
and the multiplication relation:
$$\tilde{m}_A:=\{(b_1+(1+s b_1) b_2,(1+s b_1) a_2;b_1,a_1,b_2,a_2)\,:\,a_1=(1+s b_1) (a_2-s b_2)\}
\subset\Gamma_s\times\Gamma_s\times\Gamma_s$$
Straightforward computations show that, for $s\neq 0$ the map: 
$$\Gamma_s\ni (b,a)\mapsto (-1/s,-1)(b,a)(-1/s,-1)^{-1}=(-b+(a-1)/s, a)\in\Gamma_s$$
is an isomorphism of $\Gamma_s\rra C_s$ and $\Gamma_s\rra A$. 
Groupoid structures on $\Gamma_0=G$ are given by the double Lie group $(G;A,B)$; for $s\neq 0$ 
the map:
$\Gamma_s\ni(b,a)\mapsto (0,s)(b,a)(0,s)^{-1}=(s b, a)\in \Gamma_1$ gives  isomorphisms of both groupoid structures, 
so it is enough to consider $s=1$.

Let us now denote $C:=C_1$ and $\Gamma:=\Gamma_1$. On $\Gamma$ there are two (isomorphic) groupoid structures:\\ 
$\Gamma_C: \Gamma\rightrightarrows C$ with structure
\begin{eqnarray} \tilde{c}_L(b,a)=(b,1+b)\,,\,\tilde{c}_R(b,a)=(\frac{b}{a-b},\frac{a}{a-b})\,,\,
\tilde{s}_C(b,a):=(\frac{b}{a- b},\frac{1+  b}{a-  b})\nonumber\\
\tilde{m}_C=\{(b_1, \frac{a_1 a_2}{1+ b_2};b_1,a_1,b_2,a_2)\,:\,b_1=b_2 (a_1- b_1)\}\end{eqnarray}
$\Gamma_A: \Gamma\rightrightarrows A$ with structure
\begin{eqnarray}
a_L(b,a)=(0,a-b)\,,\,a_R(b,a)=(0,\frac{a}{1+b})\,,\,s_A(b,a)=(\frac{-b}{1+b},\frac{a-b}{1+b}),\nonumber\\
m_A=\{(b_1+(1+b_1) b_2,(1+ b_1) a_2;b_1,a_1,b_2,a_2)\,:\,a_1=(1+b_1) (a_2-b_2)\}\end{eqnarray}
Since $(G;B,C)$ is a double Lie group there are two  groupoid structures on $G$:\\
$G_B: G\rra B$:
\begin{eqnarray}
b_L(b,a):=(b-a+1,1)\,, \,\, b_R(b,a):=(\frac{b-a+1}{a},1) \,,\,\,s_B(b,a):=(\frac{2 - 2 a +b}{a},\frac1a)  \nonumber\\
m_B:=\{(b_1+a_1(a_2-1), a_1 a_2;b_1,a_1,b_2,a_2): 1+ b_1 -a_1=a_1(1 + b_2 -a_2)\}\end{eqnarray}
$G_C: G\rra C$:
\begin{eqnarray}
c_L(b,a):=(a-1,a)=c_R(b,a):= \,,\,\,s_C(b,a):=(2-2 a-b,a)  \nonumber\\
m_C:=\{(b_1+b_2+1-a, a; b_1,a,b_2,a)\}\end{eqnarray}
Let us  denote $\delta_0:=m_C^T$;  this is a coassociative morphism: $G_B\rel G_B\times G_B$.

To summarize, we arrive at the following situation: there is a Lie group $G$ and three closed subgroups $A,B,C$ satisfying conditons
$G=BC\,,\,A\cap C=B\cap C=\{e\}$. 
This is investigated in the next section.
%%%%%%%%%%%%%%%%%%%%%%%%%%%%%%%%%%%%%%%%%%%%%%%%%%%%%%%%%%%%%%%%%%%%%%%%%%%%%%%%%%%%%%%%%%%%%%%%%%%
%%%%%%%%%%%%%%%%%%%%%%%%%%%%%%%%%%%%%%%%%%%%%%%%%%%%%%%%%%%%%%%%%%%%%%%%%%%%%%%%%%%%%%%%%%%%%%%%%%%%
%%%%%%%%%%%%%%%%%%%%%%%%%%%%%%%%%%%%%%%%%%%%%%%%%%%%%%%%%%%%%%%%%%%%%%%%%%%%%%%%%%%%%%%%%%%%%%%%%
%\input{rozdzial3}
\section{The twist}\label{stwist}
Let $G$ be a group and $A,B,C\subset G$ subgroups satisfying conditions: 
\notka{warunki}\begin{equation}\label{warunki} 
B\cap C=\{e\}=A\cap C\,,\,B C=G.
\end{equation}
\notka{warunki}
 As described above, in this situation, there is the groupoid $G_B$, and the (coassociative) morphism
$\delta_0:G_B\rel G_B\times G_B$; explicitely the graph of $\delta_0$ is equal to: 
\notka{delta0}\begin{equation}\label{delta0} \delta_0=\{(b_1 c, c b_2;b_1 c b_2)
  \,:\, b_1,b_2\in B\,,\,c\in C\}
\end{equation}
%%%%
Let us note that  $c_L(g),b_L(g)$ and $c_R(g),b_R(g)$ determine $g$ uniquely:
$$g=b_R(b_R(g) c_R(g)^{-1}) c_R(g)=c_L(g) b_L(c_L(g)^{-1} b_L(g))$$
Using the lemma \ref{embed} we see that this is a transformation groupoid $(C\backslash G) \times C$ and the isomorphism is
$$(C\backslash G) \times C\ni ([g],c)\mapsto b_R(g)c\in G$$
%%%%%%%%
Let $\Gamma:=A C\cap C A$ and  consider on $\Gamma$ the groupoid structure $\Gamma_A$ over $A$ described above, together with a relation
$$\tilde{m}_C^T:=\{(a_1 c_1, c_1 a_2 ; a_1 c_1 a_2) : a_1 c_1, c_1 a_2\in\Gamma\}\subset \Gamma\times\Gamma\times\Gamma.$$
The corresponding projections will be denoted by $\tilde{c}_L, \tilde{c}_R $ and 
$a_R, a_L$. Again, by the lemma \ref{embed} we identify the groupoid $\Gamma_A$ with the restriction 
of $(C\backslash G) \times C$ and then with the restriction
of $G_B$  to the set $B':=B\cap CA$, i.e. with $b_L^{-1}(A')\cap b_R^{-1}(A')$. 
This restriction will be denoted by $\Gamma_{B'}$ (instead of more adequate but rather inconvenient $G_B|_{B'}$). 
This isomorphism is given by:
$\Gamma_A\ni a c\mapsto b_R(a) c\in G_B$. 

The image of $\tilde{m}_C^T$ is equal to:
\begin{equation}\label{mct}
\{(b_R(a_1) c_1, b_R(a_2) c_2; b_R(a_1 a_2) c_2): a_1 c_1=\tilde{c}_1 \tilde{a}_1, c_1\tilde{a}_2=a_2 c_2\}
\end{equation}
%%%%%%%%%%%%%%%%%%%%%%%%%%%%%%%%%%%%%%%%%%%
Let us now define  a basic object of this section
\notka{twist}\begin{equation}\label{twist}
T:=\{(g,b):c_R(g)b\in A\}=\{(b_1 \tilde{c}_L(b_2)^{-1}, b_2) : b_1\in B,b_2\in B'\}\subset G_B\times G_B.
\end{equation}
Using the definition (\ref{delta0}) of $\delta_0$   one easily computes images of $T$ by 
relations $id\times\delta_0$ and $\delta_0\times id$:
\notka{twist-1}\begin{equation}\label{twist-1}
(id\times \delta_0)T=\{(g_1,b_2,b_3): b_2 b_3\in B'\,,\,c_R(g_1)=\tilde{c}_L(b_2 b_3)^{-1}\}
\end{equation}\vsm\vsm\vsm
\notka{twist-2}\begin{equation}\label{twist-2}
(\delta_0\times id)T=\{(g_1,g_2,b_3): c_R(g_1)=c_L(g_2)\,,\,b_3\in B'\,,\,c_R(g_2)=\tilde{c}_L(b_3)^{-1}\}
\end{equation}
Let us also denote $T_{12}:=T\times B\subset G_B\times G_B\times G_B$ and 
$T_{23}:=B\times T\subset G_B\times G_B\times G_B$. \\
%%%%%%%%%%%%%%
Main properties of $T$ are listed in the following
\begin{prop} \label{alg-twist}\notka{alg-twist}\begin{enumerate}
\item $T$ is a section of left and right projections (in $G_B\times G_B$) over the set $B\times B'$ and a 
bisection of $G_B\times \Gamma_{B'}$;
\item $(id\times \delta_0) T$ is a section of left and right projections 
(in $G_B\times G_B\times G_B$) over the set 
$B\times \delta_0(B')=\{(b_1,b_2,b_3) : b_2 b_3\in B'\}$;
\item $(\delta_0\times id) T$ is a section of left and right projections  over the set 
$B\times B \times B'$;
\item   $T_{23}(id\times \delta_0)T= T_{12}(\delta_0\times id) T$ (equality of sets in 
$G_B\times G_B\times G_B$), 
moreover this set is a section of the right projection over $B\times (\delta_0(B')\cap(B\times B'))$ and 
the left projection over
$B\times B'\times B'$.
\end{enumerate}
\end{prop}
{\em  Proof:} 1) Since the ``right leg'' of $T$ is  $B'$ it is clear that if $T$ is a section of left 
and right projection 
over $B\times B'$ then it is a bisection of $G_B\times  \Gamma_{B'}$. The formula (\ref{twist}) implies  that 
$T$ is a section of $b_L\times b_L$ over   $B\times B'$. But since 
$g\in G$ is determined by $b_R(g)$ and $c_R(g)$,  it is enough to show that 
for $(b_1,b_2)\in B\times B'$ there exists $g_1$ with $b_R(g_1)=b_1$ and 
$(g_1, b_2)\in T$.  Let $b_2=:c_0 a_0$ then for $g_1:= s_B(b_1 c_0)$ we have $b_R(g_1)=b_1$ 
and $c_R(g_1)=c_0^{-1}$ so   
$(g_1,b_2)$ is an element of  $T$.\\
2) From (\ref{twist-1}) it follows that $(id\times \delta_0) T$ is a section of the left 
projection over $B\times \delta_0(B')$.
As in point 1), for $(b_1,b_2,b_3)\in B\times \delta_0(B')$ with $b_2 b_3=:c_0 a_0$ one has 
$b_R(s_B(b_1 c_0))=b_1$ and 
$(s_B(b_1 c_0),b_2,b_3)\in (id\times \delta_0) T$.\\
3) Follows directly from (\ref{twist-2}).\\
4) By direct computation one gets:
$$T_{23}(id\times \delta_0)T=\{(g_1,g_2,b_3): b_3, b_R(g_2) b_3 \in B'\,,\,c_R(g_2)=\tilde{c}_L(b_3)^{-1}\,,\,
c_R(g_1)=\tilde{c}_L(b_R(g_2)b_3)^{-1}\}$$
and the same result for $T_{12}(\delta_0\times id) T$. This formula implies that for any 
$(b_1,b_2,b_3)\in  B\times (\delta_0(B')\cap(B\times B'))$ there exists exactly one pair $(g_1,g_2)$ such that
$b_1=b_R(g_1), b_2=b_R(g_2)$ and $(g_1,g_2, b_3)\in T_{23}(id\times \delta_0)T$. 
It remains to prove that it is a section of the left projection over $B\times B'\times B'$. It is clear that 
for a given $(b_1,b_2,b_3)$  there exists at most one pair $(g_1,g_2)$ such that $b_1=b_L(g_1), b_2=b_L(g_2)$ and 
$(g_1,g_2,b_3)\in T_{23}(id\times \delta_0)T$. The rest follows from the following observation:
for  $b_3\in B'\,,\,b_3=:c_0 a_0$, there is an equivalence 
$$b\in B'\iff b_R(bc_0^{-1})b_3\in B'$$

Indeed if $b=:\tilde{c}_0\tilde{a}_0$ then 
$$b_R(b c_0^{-1})b_3=c_L(b c_0^{-1})^{-1} b c_0^{-1} b_3=c_L(b c_0^{-1})^{-1} b a_0 =
c_L(b c_0^{-1})^{-1} \tilde{c}_0\tilde{a}_0 a_0\in C A$$
On the other direction,  
$$b_R(bc_0^{-1})b_3\in C A\Rightarrow bc_0^{-1}b_3\in C A\Rightarrow b a_0\in C A\Rightarrow b\in C A$$
\dow

 If a subset $S$ of a groupoid $\Gamma$ is a section of the right projection (over $e_R(S)$) then every element of 
$e_L^{-1}(e_R(S))$ can be multiplied from the left by unique  element of $S$, therefore $S$ defines a mapping 
$e_L^{-1}(e_R(S))\rightarrow e_L^{-1}(e_L(s))$ -- the left multiplication by $S$. If $S$ is also a section of the left 
projection, this mapping is a bijection.
Because of the prop \ref{alg-twist} the left multiplication by $T$, which we denote by the same symbol, is a 
bijection of
$G_B \times b_L^{-1}(B')$; also the left multiplication by $(\delta_0\times id)T$, 
which will be denoted by $T_1$, is a bijection of $G_B\times G_B \times b_L^{-1}(B')$, 
and the  left multiplication by $(id\times \delta_0)T$, denoted by $T_2$, is a bijection of 
$G_B\times (b_L\times b_L)^{-1}(\delta_0(B'))$. These mappings are given by:
\notka{T}\begin{eqnarray}\label{T}
T: (b_1 c_1,b_2 c_2)\mapsto (s_B(b_1 \tilde{c}_L(b_2)) c_1, b_2
c_2) =(c_L(b_1 \tilde{c}_L(b_2))^{-1} b_1 c_1,b_2 c_2)=\\
\nonumber =(b_R(b_1 \tilde{c}_L(b_2))\tilde{c}_L(b_2)^{-1}c_1, b_2 c_2)
\end{eqnarray}
$$T^{-1}: (b_1 c_1,b_2 c_2)\mapsto (s_B(b_1 \tilde{c}_L(b_2)^{-1})
c_1, b_2 c_2)$$\vspace{-3ex}
\notka{T1}\begin{eqnarray}\label{T1}
T_1:(b_1 c_1,b_2 c_2,b_3 c_3)\mapsto (s_B(b_1 c_L(b_2\tilde{c}_L(b_3))) c_1, s_B(b_2\tilde{c}_L(b_3))c_2, b_3 c_3) =\\
\nonumber=( c_L(b_1 b_2 \tilde{c}_L(b_3))^{-1}b_1 c_1 , 
c_L(b_2 \tilde{c}_L(b_3))^{-1}b_2 c_2,b_3 c_3 )
\end{eqnarray}
$$T_1^{-1}:(b_1 c_1,b_2 c_2,b_3 c_3)\mapsto (s_B(b_1
c_L(b_2\tilde{c}_L(b_3)^{-1})) c_1, s_B(b_2\tilde{c}_L(b_3)^{-1})c_2, b_3 c_3)$$\vspace{-3ex}
\notka{T2}\begin{eqnarray}\label{T2}
T_2:(b_1 c_1,b_2 c_2,b_3 c_3)\mapsto (s_B(b_1\tilde{c}_L(b_2 b_3))
c_1, b_2 c_2, b_3 c_3)=\\\nonumber
=( c_L(b_1\tilde{c}_L(b_2 b_3))^{-1} b_1 c_1, b_2 c_2, b_3 c_3)
\end{eqnarray}
$$T_2^{-1}:(b_1 c_1,b_2 c_2,b_3 c_3)\mapsto (s_B(b_1\tilde{c}_L(b_2 b_3)^{-1})
c_1, b_2 c_2, b_3 c_3)$$
%%%%%
The left multiplication by the composition $T_{12}T_1=T_{23}T_2$ is a 
bijection from $G_B\times (b_L\times b_L)^{-1}(\delta_0(B')\cap(B\times B'))$ to
$G_B\times  (b_L\times b_L)^{-1}( B'\times B')$ and is given by:
$$T_{12}T_1: (b_1c_1,b_2c_2,b_3c_3)\mapsto (c_L(b_1\tilde{c}_L(b_2 b_3))^{-1}b_1 c_1, 
c_L(b_2\tilde{c}_L(b_3))^{-1}b_2 c_2,b_3c_3)$$

Let $Ad_{T}:G_B\times G_B\rel G_B\times G_B$ be a relation defined by:
$$(g_1,g_2;g_3,g_4)\in Ad_{T}\iff \exists t_1, t_2\in T : (g_1,g_2)=t_1(g_3,g_4)(s_B\times s_B)(t_2).$$
Using  the definition of $T$ (\ref{twist}) one gets: 
$$Ad_T=\{(b_R(b_3\tilde{c}_L(b_4))\tilde{c}_L(b_4)^{-1} c_3 \tilde{c}_L(b_R(b_4 c_4)),b_4 c_4; 
b_3 c_3, b_4 c_4): b_4,b_R(b_4 c_4)\in B'\}.$$  
Finally, let us define  the relation $\delta:=Ad_T\cdot \delta_0: G_B\rel G_B\times G_B$.
\notka{delta}\begin{equation}\label{delta} 
\delta=\{(b_R(b_3 b_2^{-1} \tilde{c}_L(b_2))\tilde{c}_L(b_2)^{-1} c_L(b_2 c_2) \tilde{c}_L(b_R(b_2 c_2)), 
b_2 c_2; b_3 c_2) : b_2,b_R(b_2 c_2)\in B'\}
\end{equation}
The next lemma explains the relation between $\delta$ and $\tilde{m}_C^T$:
\begin{lem} $\delta$ is an extension of $\tilde{m}_C^T$ i.e. $\tilde{m}_C^T\subset \delta$
\end{lem}
{\em Proof:} Recall that (\ref{mct}):
$$\tilde{m}_C^T=\{(b_R(a_1) c_1, b_R(a_2) c_2; b_R(a_1 a_2) c_2): a_1 c_1=\tilde{c}_1 \tilde{a}_1, a_2 c_2=c_1\tilde{a}_2\}$$
By the lemma \ref{embed} $b_R(a_2)$ and $b_R(b_R(a_2) c_2)$ are in $B'$. So it remains to prove that:
$$b_R\left[b_R(a_1 a_2) b_R(a_2)^{-1} \tilde{c}_L(b_R(a_2))\right]\tilde{c}_L(b_R(a_2))^{-1} c_L(b_R(a_2) a_2) \tilde{c}_L(b_R(b_R(a_2) c_2))=
b_R(a_1) c_1$$
Let $a_2=c_0 b_0$ then $b_R(a_2)=b_0$ and  $\tilde{c}_L(b_R(a_2))=c_0^{-1}$. We have:
$$b_R\left[b_R(a_1 a_2) b_R(a_2)^{-1} \tilde{c}_L(b_R(a_2))\right]=b_R(a_1 a_2 b_0^{-1} c_0^{-1})=b_R(a_1)$$
and
$$b_R(a_1)\tilde{c}_L(b_R(a_2))^{-1} c_L(b_R(a_2) c_2) \tilde{c}_L(b_R(b_R(a_2) c_2))=b_R(a_1) c_0 c_L(b_0 c_2) \tilde{c}_L(b_R(b_R(a_2) c_2))=$$
$$=b_R(a_1) c_L(c_0 b_0 c_2) \tilde{c}_L(b_R(c_0 b_0 c_2))=b_R(c_1)\tilde{c}_L( c_L(c_0 b_0 c_2)b_R(c_0 b_0 c_2))=$$
$$=b_R(a_1)\tilde{c}_L( a_2 c_2))=b_R(a_1) c_1$$
\dowl

Everything above was purely algebraic. Now we add some differential conditions, 
that enable us to use $T$ to twist the comultiplication on $C^*_r(G_B)$. 

{\em The following  are standing assumptions for the rest of this section}
\begin{assumpt}\label{zal}\notka{zal}
\begin{enumerate}
\item $G$ is a Lie group and $A,B,C$ are closed Lie subgroups such that \\
$B\cap C=\{e\}=A\cap C\,,\,B C=G.$
\item The set $\Gamma:=C A\cap AC$ is open and dense in $G$.
\item Let $U:=b_L^{-1}(B')$ and $\sA(U)$ be the linear space of elements from $\sA(G_B)$ supported in $U$.
We assume that $\sA(U)$ is dense in $C^*_r(G_B)$.
\item For a compact set $K_C\subset C$, open $V\subset B$ and $(b_1, b_2)\in B\times B'$ let us define a set
$Z(b_1,b_2,K_C;V):=K_C\cap\{c\in C: b_R(b_1 c) b_2\in V\}$ and a function:
$$B\times B' \ni(b_1,b_2)\mapsto \mu(b_1,b_2,K_C;V):=\int_{Z(b_1,b_2,K_C;V)} d_l c.$$
For compact sets $K_1\subset B$ and $K_2\subset B'$ let $\mu(K_1,K_2,K_C;V):=\sup\{\mu(b_1,b_2,K_C;V): b_1\in K_1\,,\,b_2\in K_2\}$
We assume that \begin{equation}
\forall\,\epsilon>0\,\exists\, V-{\rm a\, neighbourhood\,of\,}B\setminus B'{\rm\, in\, B}\,:\, \mu(K_1,K_2,K_B;V)\leq\epsilon
\end{equation}
\end{enumerate}
\end{assumpt}
These technical assumptions are sufficient to prove all we need and are satisfied in examples we are 
intersted in (probably they can be weakened and some of them imply others). 
The strategy is to interpret  all in  $C^*_r(G_B)$ where there is  a well defined comultiplication.
The first simple result is:

\begin{lem}
$C^*_r(\Gamma_{B'})=C^*_r(G_B)$
\end{lem}
{\em Proof}: Let $U$ and $\sA(U)$ be as in assumptions \ref{zal}. Then $J:=\sA(U)$ is a right ideal in 
$\sA(G_B)$. Then $J^*$ is a left  ideal which is also dense in $C^*_r(G_B)$. 
The subalgebra $J J^*$ is contained in $\sA(\Gamma_{B'})$. Every positive element in $C^*_r(G_B)$ can be 
approximated by elements from $J J^*$, so the same is true for every element in  $C^*_r(G_B)$.\\\dowl

It is straightforward to check that assumptions \ref{zal} guarantee
that all sets appearing in prop \ref{alg-twist} are submanifolds and
the corresponding mappings are diffeomorphisms. Since $T$ is a
bisection over the open set $B\times B'$ it defines (by a push-forward)
mapping of $\sA(G_B\times U)$ which will be denoted again by $T$. The
prop \ref{alg-twist} suggests that it can be used to twist the
comultiplication $\Delta_0$ on $C^*_r(G_B)$. This is the content of
the next proposition.

%%%%%%%%%%%%%
\begin{prop}\label{Delta} \notka{Delta} a) The mapping
  $T:\sA(G_B\times U)\rightarrow\sA(G_B\times U)$
extends to a unitary $\widehat{\sT}\in M(C^*_r(G_B)\mt C^*_r(G_B))$ which satisfies:
$$(\widehat{\sT}\mt I)(\Delta_0\mt id)\widehat{\sT}=(I \mt \widehat{\sT})(id\mt \Delta_0)\widehat{\sT}$$
b) Because of a), the formula $\Delta(a):=\widehat{\sT}\Delta_0(a)\widehat{\sT}^{-1}$ defines a coassociative morphism. 
For this morphism:
$$cls\{\Delta(a)(I\mt c)\,,\,a,c\in C^*_r(G_B)\}=cls\{\Delta(a)(c\mt I)\,,\,a,c\in C^*_r(G_B)\}=
C^*_r(G_B)\mt C^*_r(G_B)$$
("cls" stands for "closed linear span").
\end{prop}

The rest of this section is dedicated to the proof of this proposition.

\noindent{\em Proof of the statement a):}\\
Mappings $T, T_1, T_2, T_{12}, T_{23}, T_{12}T_1=T_{23}T_2 $ are originally defined on linear subspaces of 
$\sA(G_B\times G_B)$ and $\sA(G_B\times G_B\times G_B)$. So the first task is to extend them to multipliers.
To this end we will use the following simple lemma:
\begin{lem} \label{mnoznik} \notka{mnoznik} Let $A$ be a $C^*$-algebra  and 
$V_1,V_2$ be dense linear subspaces of  $A$. Suppose that a linear bijection
$T:V_1\rightarrow V_2$ satisfies $ \om_2^*(T \om_1)=(T^{-1}\om_2)^* \om_1$ for any $\om_2\in V_2,\,\om_1\in V_1$. 
Then $T$ extends to a unitary multiplier of $A$.
\dowl
\end{lem}
We need also the lemma which is a straightforward generalization of the similar fact for a 
bisection which was proven in \cite{DG}.
\begin{lem} \label{T-conjugate} \notka{T-conjugate}
Let $\Gamma$ be a differential groupoid and $C\subset\Gamma$ be a submanifold such that  
$\el|_C:C\rightarrow C_l$ and $\er|_C:C\rightarrow C_r$ are 
diffeomorphisms onto open sets $C_r,C_l\subset E$. Then for $\om_1\in\sA(\el^{-1}(C_r))\,,\,\om_2\in\sA(\el^{-1}(C_l))$ we have:
$\displaystyle \om_2^*(C\om_1)=(s(C)\om_2)^*\om_1$ and $\displaystyle (C\om_1)^*(C\om_1)=\om_1^*\om_1$.
\dowl
\end{lem}
\noindent Using the lemma we obtain equalities:
\begin{eqnarray}
\nonumber \om_2^*(T\om_1)=(T^{-1}\om_2)^*\om_1\,,\,\,\,\om_1,\om_2\in\sA(G_B\times U)\,,\\
\nonumber \om_2^*(T_1\om_1)=(T_1^{-1}\om_2)^*\om_1\,,\,\,\,\om_1,\om_2\in \sA(G_B\times G_B\times U)\,,\\
\nonumber \om_2^*(T_2\om_1)=(T_2^{-1}\om_2)^*\om_1\,,\,\,\,\om_1,\om_2\in \sA(U_2)\,,\\
\nonumber \om_2^*(T_{23}\om_1)=(T_{23}^{-1}\om_2)^*\om_1\,,\,\,\,\om_1,\om_2\in\sA(G_B\times G_B\times U)\,,\\
\nonumber \om_2^*(T_{12}\om_1)=(T_{12}^{-1}\om_2)^*\om_1\,;\,\,\om_1,\om_2\in\sA(G_B\times U \times G_B),\end{eqnarray}
where $U_2:=G_B\times (b_L\times b_L)^{-1}(\delta_0(B'))$

Now, by assumptions \ref{zal} (3)   it is clear that $\sA(G_B\times U )\,,\,
\sA(G_B\times G_B \times U)$ and $\sA(G_B\times U \times G_B)$  are dense in 
$C^*_r(G_B)\mt C^*_r(G_B)$ and $C^*_r(G_B)\mt C^*_r(G_B)\mt C^*_r(G_B)$ respectively, 
so we have (unitary) multipliers $\widehat{\sT}$, $\widehat{T_1}$, $\widehat{T_{12}}$ and $\widehat{T_{23}}$.
 
To extend $T_2$ we need density of  $\sA(U_2)$ i.e density of $\sA((b_L\times b_L)^{-1}(\delta_0(B'))$.
By the results  of \cite{DLG}, the set $\widehat{\delta_0}(\sA(G_B))\sA(G_B\times G_B)$ is linearly dense in $C^*_r(G_B)\mt C^*_r(G_B)$, 
since $\sA(U)$  is dense, the same is true for the set 
$\widehat{\delta_0}(\sA(U))\sA(G_B\times G_B )$, 
using formula (\ref{delta0-form}) one checks that this set is contained in 
$\sA((b_L\times b_L)^{-1}(\delta_0(B'))$. Therefore $T_2$ extends to 
unitary multiplier.

Now, we know that (unitary) multipliers  $\displaystyle \widehat{T_{12}}\widehat{T_1}$ and $\widehat{T_{23}}\widehat{T_2}$ restricted to 
the linear space  $\sA(G_B\times (b_L\times b_L)^{-1}(\delta_0(B')\cap(A\times B')))$ are equal and define bijection to $\sA(G_B\times U\times U)$.
Since the last set is dense in $C^*_r(G_B)\mt C^*_r(G_B)\mt C^*_r(G_B)$ they must be equal.

Now we have multipliers $\widehat{\sT}$, $\widehat{T_1}$, $\widehat{T_2}$, $\widehat{T_{12}}$ and $\widehat{T_{23}}$ together with the equation
$\displaystyle \widehat{T_{12}}\widehat{T_1}=\widehat{T_{23}}\widehat{T_2}$.
So to prove the statement a) of the proposition, it remains to show equalities: 
$$\widehat{T_{12}}=(\widehat{\sT}\mt I)\,,\,\, \widehat{T_1}=(\Delta_0\mt id)\widehat{\sT}\,,\,\,
\widehat{T_{23}}=(I \mt \widehat{\sT})\,,\,\, \widehat{T_2}=(id\mt \Delta_0)\widehat{\sT}.$$
For $\widehat{T_{23}}$ and $(I \mt \widehat{\sT})$ this is straightforward, since both multipliers agree on a dense set 
$\sA(G_B)\mt\sA(G_B\times U )$.
The same arguments work for $\widehat{T_{12}}$ and $(\widehat{\sT}\mt I)$ because they agree on a dense set 
$\sA(G_B\times U)\mt\sA(G_B)$.
For two other equalities we will use two lemmas that will be proven in the end of  the proof of the statement a) of the proposition.
%%%%%%%%%%%%%%%
\begin{lem}\label{product-mor}\notka{product-mor}
Let $\Gamma_1,\Gamma_2,\Delta_1,\Delta_2$ be differential groupoids, and 
$h_1:\Gamma_1\rel\Delta_1\,, \,h_2:\Gamma_2\rel\Delta_2$ be  morphisms. Assume that  mappings 
$\widehat{h}_1,\widehat{h}_2$ extends to morphisms $\phi_1,\phi_2$ of corresponding reduced $C^*$-algebras. 
Then $(\phi_1\mt\phi_2)(\om)=(\widehat{h_1\times h_2})(\om)$ for $\om\in\sA(\Gamma_1\times\Gamma_2)$
\end{lem}
\begin{lem} \label{id-delta-TT} \notka{id-delta-TT}For
  $\om_1\in\sA(G_B\times U)$ and 
$\om_2\in\sA(G_B\times G_B\times G_B)$:
$$T_1[\widehat{(\delta_0\times id)}(\om_1)\om_2]=\widehat{(\delta_0\times id )}(T\om_1)\om_2,$$
$$T_2[\widehat{(id\times\delta_0)}(\om_1)\om_2]=\widehat{(id\times\delta_0)}(T\om_1)\om_2.$$
\end{lem}
%{\bf sprawdzone dla $F_1$}

Now, the multiplier $(id\mt\Delta_0)\widehat{\sT}$ is defined by 
$$[(id\mt\Delta_0)\widehat{\sT}]((id\mt\Delta_0)(a)
b):=(id\mt\Delta_0)(\widehat{\sT} a)b\,,\,a\in C^*_r(G_B\times G_B),
b\in C^*_r(G_B\times G_B\times G_B).$$ 
In fact, it is enough to take $a=\om_1\in\sA(G_B\times U)$ 
and $b=\om_2\in\sA(G_B\times G_B\times G_B)$. Let us compute:
$$[(id\mt\Delta_0)\widehat{\sT}]((id\mt\Delta_0)(\om_1)\om_2)=(id\mt\Delta_0)(\widehat{\sT} \om_1)\om_2=$$
$$=\widehat{(id\times\delta_0)}(T\om_1)\om_2=
T_2[\widehat{(id\times\delta_0)}(\om_1)\om_2]=T_2[(id\mt\Delta_0)(\om_1)\om_2]$$

Therefore $(id\mt\Delta_0)\widehat{\sT}=\widehat{T_2}$. In the same way the second equality can be proven. So  
to complete the proof of the statement a) it remains to prove lemmas  \ref{product-mor}, \ref{id-delta-TT}\\
{\em Proof of lemmas:}

\noindent {\em Lemma \ref{product-mor}:} It is known \cite{DLG}, that for $\om\in\sA(\Gamma_1\times\Gamma_2)$ the mapping 
$(\widehat{h_1\times h_2})(\om)$ defines a multiplier of 
$C^*_r(\Delta_1\times\Delta_2)=C_r^*(\Delta_1)\mt C_r^*(\Delta_2)$. 
Let $\pi_{12}$ denotes the representation defined by $h_1\times h_2$ on $H:=L^2(\Delta_1\times \Delta_2)$.
For $\om'\in\sA(\Delta_1\times\Delta_2)$ we have:
$$||(\widehat{h_1\times h_2})(\om) \om'||_{C_r^*}=||\pi_{12}(\om)\pi_{id}(\om')||_{B(H)}
\leq||\pi_{12}(\om)||_{B(H)}\,||\pi_{id}(\om')||_{B(H)}=
||\pi_{12}(\om)||_{B(H)}\,||\om'||_{C_r^*}$$
Therefore we have the estimate:
$$||(\widehat{h_1\times h_2})(\om)||_{M(C_r^*)}\leq ||\pi_{12}(\om)||_{B(H)}\leq ||\om||_0$$
Now take a sequence  $\om_n\in\sA(\Gamma_1)\mt\sA(\Gamma_2)$ that converges to $\om$ in a norm $||\cdot||_0$, so also in a 
$C^*_r(\Gamma_1)\mt C^*_r(\Gamma_2)$. So $(\widehat{h_1\times h_2})(\om_n)$ converges to 
$(\widehat{h_1\times h_2})(\om)$ in $M(C^*_r(\Delta_1)\mt C^*_r(\Delta_2))$. 
On the other hand $(\widehat{h_1\times h_2})(\om_n)=(\phi_1\mt\phi_2)(\om_n)$ so it converges to $(\phi_1\mt\phi_2)(\om)$.
\\\dowl

\noindent{\em Lemma \ref{id-delta-TT}:}  To prove those formulae we need the explicit form of the action of $T$, $T_1$ and $T_2$. 
%%%%%%%%%%%%%%            FORMULA NA T   %%%%%%%%%%%%%%%%%%%%%%%%%
Let us start  with $T$. Writing 
$\om=:F\om_0$ for $F\in \sD(G_B\times U)$ we define:
$T(F \om_0)=: (TF) \om_0$. As for any bisection   $TF$ is given by \cite{DG} :
$$(TF)(T(g_1,g_2))=F(g_1,g_2)\frac{(\rho_0\mt\rho_0) (v g_1\mt w
  g_2)}{(\rho_0\mt\rho_0)(T(v g_1\mt w g_2))}\,,
\,v,w\in\Lambda^{max}(T_eC)$$
Let $(g_1, g_2)=(b_1 c_1, b_2 c_2)$ and  $T(g_1, g_2)=:(\tilde{g}_1 , g_2)$. Using
the definition of $T$ we obtain, for a curve $c_0(t)\subset C$ with
$c_0(0)=e$, 
$T(c_0(t) g_1, g_2)=(c_3(t) \tilde{g}_1, g_2)$, where 
$$c_3(t)= c_L\left[s_B(b_1\tilde{c}_L(b_2)) b_L(b_1^{-1} c_0(t) b_1)
  s_B(b_1\tilde{c}_L(b_2))^{-1}\right]^{-1}
c_L(b_1\tilde{c}_L(b_2))^{-1} c_0(t) c_L(b_1\tilde{c}_L(b_2)),$$
and $T( g_1, c_0(t) g_2)=(c_4(t) \tilde{g}_1 , c_0(t) g_2)$.\\
Identifying tangent spaces to right fibers with $T_eC$ in corresponding points,
one sees that the map we have to consider has the form 
$\left(\begin{array}{cc} M_1 & M_2\\0 & I\end{array}\right)$, so its
action on densities is determined by $M_1$ i.e. (derivative of)
$c_0(t)\mapsto c_3(t)$ and this is given by:
$$M_1: T_eC\ni\dot{c}\mapsto - P_C\cdot Ad(s_B(b_1\tilde{c}_L(b_2)))\cdot
P_B\cdot Ad(b_1)^{-1} \dot{c}+
Ad(c_L(b_1\tilde{c}_L(b_2)))^{-1}\dot{c},$$
where $P_C, P_B$ are projections in a Lie algebra $\gotg$ corresponding to the decomposition $\gotg=\gotb\oplus\gotc$.
Straightforward computation gives:
$M_1=Ad(c_L(b_1\tilde{c}_L(b_2)))^{-1}|{\gotc}.$  Finally using the
definition of $\rho_0$ and modular functions (\ref{modular-functions}) we obtain
\notka{T-form}\begin{eqnarray}\label{T-form}
(TF)(g_1,g_2)=F(T^{-1}(g_1,g_2)) j_C(\tilde{c}_L(b_L(g_2)))^{\frac12},\\
(T^{-1}F)(g_1,g_2)=F(T(g_1,g_2)) j_C(\tilde{c}_L(b_L(g_2)))^{-\frac12}\nonumber
\end{eqnarray}

%%%%%%%%%%%%%%%   FORMULA NA T_2          %%%%%%%%%%%%%%%%%%%%%%%%%%%%%

$T_2$ defines a map, again denoted by the same letter,
of $\sA(G_B\times (b_L\times b_L)^{-1}(\delta_0(B')))$. In the same way
as for $T$ we have to compute the action of $T_2$ on
$\rho_0\mt\rho_0\mt\rho_0$.
So let $(g_1,g_2,g_3)=:(b_1 c_1, b_2 c_2, b_3 c_3)$ be in the domain
of $T_2$ and $T_2(g_1, g_2, g_3)=(s_B(b_1 \tilde{c}_L(b_2 b_3)) c_1, g_2
, g_3)=:(\tilde{g}_1,g_2,g_3)$. For a curve $c_0(t)\subset C$ with
$c_0(0)=e$,  we obtain:
$T_2(c_0(t) g_1, g_2, g_3)=(c_1(t)\tilde{g}_1, g_2, g_3)$ and 
$$c_1(t)=c_R\left[c_L(b_1 \tilde{c}_L(b_2 b_3))^{-1} b_1 c_R(b_1^{-1}
  c_0(t) b_1) b_1^{-1}  c_L(b_1 \tilde{c}_L(b_2 b_3))\right]$$
Reasoning in the same way as for $T$, the action on densities is
determined  by the map:
$$T_eC\ni\dot{c}\mapsto  P_C\cdot Ad(s_B(b_1\tilde{c}_L(b_2 b_3)))\cdot
P_C\cdot Ad(b_1)^{-1} \dot{c}$$
Using the definition of modular functions (\ref{modular-functions}) one gets that the absolute value of the determinanant
of this map is equal to: $j_C(b_R(c_1 \tilde{c}_L(b_2 b_3)))\left[j_C(b_1) j_C(\tilde{c}_L(b_2 b_3))\right]^{-1}$.
In this way, we get:
\notka{T2-form}\begin{eqnarray}\label{T2-form}
(T_2F)(g_1, g_2, g_3)=F(T_2^{-1}(g_1, g_2, g_3))
j_C(\tilde{c}_L(b_2 b_3))^{\frac12}\\
(T_2^{-1}F)(g_1, g_2, g_3)=F(T_2(g_1, g_2, g_3)) j_C(\tilde{c}_L(b_2 b_3))^{-\frac12}\nonumber
\end{eqnarray}

%{\bf Formuly na  $T_1$}

Finally, similar computations as for $T_2$ give for $T_1$:
\notka{T1-form}\begin{eqnarray}\label{T1-form}
(T_1 F) (g_1, g_2, g_3)=F(T_1^{-1}(g_1, g_2, g_3)) j_C(\tilde{c}_L(b_3))^{\frac12} j_C(c_L(b_2
\tilde{c}_L(b_3)^{-1}))^{-\frac12} \\
(T_1^{-1} F) (g_1, g_2, g_3)=F(T_1(g_1, g_2, g_3)) j_C(\tilde{c}_L(b_3))^{-\frac12} j_C(c_L(b_2
\tilde{c}_L(b_3)))^{\frac12} \nonumber
\end{eqnarray}

Now we can prove formulae in the lemma. Let us compute the left hand side of the first  equality:

$$T_1[(\widehat{\delta_0\times id})(F_1)F_2](b_1 c_1,b_2 c_2,b_3 c_3)=$$
$$=[(\widehat{\delta_0\times id})(F_1)F_2]
\left( s_B(b_1 c_L(b_2\tilde{c}_L(b_3)^{-1})) c_1, 
s_B(b_2\tilde{c}_L(b_3)^{-1}) c_2, b_3 c_3 \right)
\frac{j_C(\tilde{c}_L(b_3))^{\frac12}}{ j_C(c_L(b_2 \tilde{c}_L(b_3)^{-1}))^{\frac12}} =$$
$$=\frac{j_C(\tilde{c}_L(b_3))^{\frac12}}{j_C(c_L(b_2 \tilde{c}_L(b_3)^{-1}))^{\frac12}}
\int_{C\times C} d_lc' d_lc'' j_C(c_L(b_R(b_2 \tilde{c}_L(b_3)^{-1}) c'))^{-\frac12} 
F_1(b_R(b_1 b_2 \tilde{c}_L(b_3)^{-1})c', b_3 c'')\times$$ 
$$ \times F_2\left(c_L(b_1 b_2 \tilde{c}_L(b_3) c')^{-1} b_1 c_1, 
b_R(b_2  \tilde{c}_L(b_3)^{-1} c') c'^{-1} \tilde{c}_L(b_3) c_2, c_L(b_3 c'')^{-1} b_3 c_3\right)
=$$
$$=j_C(\tilde{c}_L(b_3))^{\frac12} \int_{C\times C} d_lc' d_lc''
j_C(c_L(b_2 \tilde{c}_L(b_3)^{-1}c'))^{-\frac12}
F_1(b_R(b_1 b_2 \tilde{c}_L(b_3)^{-1})c', b_3 c'')\times$$ 
$$ \times F_2\left(c_L(b_1 b_2 \tilde{c}_L(b_3) c')^{-1} b_1 c_1, 
b_R(b_2  \tilde{c}_L(b_3)^{-1} c') c'^{-1} \tilde{c}_L(b_3) c_2, c_L(b_3 c'')^{-1} b_3 c_3\right)$$
The first equality uses formula (\ref{T1-form}) for $T_1$ , the second  one the formula (\ref{delta0-id-form}) 
and the  third one  the equality:
$$j_C(c_L(b_2 \tilde{c}_L(b_3)^{-1}))j_C(c_L(b_R(b_2 \tilde{c}_L(b_3)^{-1}) c'))= j_C(c_L(b_2 \tilde{c}_L(b_3)^{-1} c'))$$
And the right hand side:
$$[(\widehat{\delta_0\times id})(T F_1)F_2](b_1 c_1,b_2 c_2,b_3 c_3)=$$
$$=\int_{C\times C} d_lc' d_lc''(TF_1)(b_1 b_2 c', b_3 c'') 
F_2(c_L(b_1 b_2 c')^{-1} b_1 c_1, b_R(b_2 c') c'^{-1} c_2, c_L(b_3 c'')^{-1} b_3 c_3) 
j_C(c_L(b_2 c'))^{-\frac12}=$$
$$=\int_{C\times C} d_lc' d_lc'' j_C(c_L(b_2 c'))^{-\frac12}j_B(\tilde{c}_L(b_3))^{\frac12} 
F_1(b_R(b_1 b_2 \tilde{c}_L(b_3)^{-1}) \tilde{c}_L(b_3) c', b_3 c'') \times$$
$$\times F_2(c_L(b_1 b_2 c')^{-1} b_1 c_1, b_R(b_2 c') c'^{-1} c_2, c_L(b_3 c'')^{-1} b_3 c_3) 
$$
After the change of variables  $(c', c'')\mapsto ( \tilde{c}_L(b_3) c', c'')$ in the last integral, 
we get the equality.

To prove the second formula we compute using formulae (\ref{T2-form}) for $T_2$ and (\ref{id-delta0-form}):
$$T_2[(\widehat{id \times \delta_0})(F_1)F_2](b_1 c_1,b_2 c_2,b_3 c_3)=$$
$$=j_C(\tilde{c}_L(b_2 b_3))^{\frac12}[(\widehat{id \times \delta_0})
(F_1)F_2](s_B(b_1 \tilde{c}_L(b_2 b_3)^{-1}) c_1, b_2 c_2, b_3 c_3)=$$
$$=j_C(\tilde{c}_L(b_2 b_3))^{\frac12} \int_{C\times C} d_lc' d_l c''\,j_C(c_L(b_3 c''))^{-\frac12} 
 F_1(b_R(b_1 \tilde{c}_L(b_2 b_3)^{-1})c', b_2 b_3 c'')  \times$$ 
$$\times 
F_2\left(c_L(b_R(b_1 \tilde{c}_L(b_2 b_3)^{-1})c')^{-1} b_R(b_1 \tilde{c}_L(b_2 b_3)^{-1}) 
\tilde{c}_L(b_2 b_3) c_1, c_L(b_2 b_3 c'')^{-1} b_2 c_2, b_R(b_3 c'')c''^{-1} c_3\right) $$
The first argument of $F_2$ reads:
$$c_L(b_R(b_1 \tilde{c}_L(b_2 b_3)^{-1})c')^{-1} b_R(b_1 \tilde{c}_L(b_2 b_3)^{-1}) 
\tilde{c}_L(b_2 b_3) c_1=$$
$$=c_L(b_R(b_1 \tilde{c}_L(b_2 b_3)^{-1})c')^{-1} b_R(b_1 \tilde{c}_L(b_2 b_3)^{-1}) 
\tilde{c}_L(b_2 b_3) c_1=$$
$$=c_L(b_R(b_1 \tilde{c}_L(b_2 b_3)^{-1})c')^{-1}  c_L(b_1 \tilde{c}_L(b_2 b_3)^{-1}) b_1 c_1=
c_L(b_1 \tilde{c}_L(b_2 b_3)^{-1} c')^{-1} b_1 c_1$$
%%%%%%%%%
Now, the right hand side:
$$[(\widehat{id \times \delta_0})(T F_1)F_2](b_1 c_1,b_2 c_2,b_3 c_3)=
\int_{C\times C} d_lc' d_l c'' \,j_C(c_L(b_3 c''))^{-\frac12}\times $$
$$\times (TF_1)(b_1 c', b_2 b_3 c'') 
F_2(c_L(b_1 c')^{-1} b_1 c_1, c_L(b_2 b_3 c'')^{-1} b_2 c_2, b_R(b_3
c'')^{-1} c''^{-1} c_3)= $$
$$=\int_{C\times C} d_lc' d_l c'' \,j_C(c_L(b_3 c''))^{-\frac12} j_C(\tilde{c}_L(b_2 b_3))^{\frac12}\times$$
$$\times F_1(c_L(b_1 \tilde{c}_L(b_2 b_3)^{-1})^{-1} b_1 c', b_2 b_3 c'') F_2(c_L(b_1 c')^{-1} b_1 c_1, 
c_L(b_2 b_3 c'')^{-1} b_2 c_2, b_R(b_3 c'')c''^{-1} c_3)$$
Now, after the change the variables in this integral $(c',c'')\mapsto (\tilde{c}_L(b_2 b_3)c', c'')$ 
we get the equality.
\\\dowl

\noindent
In this way the proof of the statement a) of the prop.
\ref{Delta} is completed.\vspace{1ex}

Now we pass to statement b) i.e. density conditions for $\Delta$. 
Since $\Delta(a)(I\mt b)=\widehat{\sT}\Delta_0(a)\widehat{\sT}^{-1}(I\mt b)$ and $\widehat{\sT}$ is
a unitary multiplier, it is enough to prove that $\Delta_0(a)\widehat{\sT}^{-1}(I\mt b)$ is linearly dense in 
$C_r^*(G_B)\mt C_r^*(G_B)$; the same is true for $\Delta_0(a)\widehat{\sT}^{-1}(b\mt I)$. 
As in \cite{DLG} we give explicit formulae for
$\Delta_0(a)\widehat{\sT}^{-1}(I\mt b)$ and $\Delta_0(a)\widehat{\sT}^{-1}(b\mt I)$ for $a,b\in\sA(G_B)$ supported in some subsets. 
Then, by continuity we obtain inclusions $\Delta_0(a)\widehat{\sT}^{-1}(I\mt b)\subset C_r^*(G_B)\mt C_r^*(G_B)$ and 
$\Delta_0(a)\widehat{\sT}^{-1}(b\mt I)\subset C_r^*(G_B)\mt C_r^*(G_B)$ for all $a,b\in C_r^*(G_B)$. Finally we prove density conditions.

To get the  formula,  using change of variables, we  write 
$(\widehat{\delta_0}(f_1)T^{-1}(f_3\mt f_2*f_4))$ as $K_1(f_1,f_2)*(f_3 \mt f_4)$ for some function 
$K_1(f_1,f_2)\in \sD(G_B\times G_B)$ at least for $f_1,f_2,f_3,f_4$ in some 
large enough subspace of functions. Let us define sets and mappings:
%%%%
\notka{dfi1}\begin{equation}\label{dfi1}
D(\Phi_1):=\{(b_1 c_1 , b_2 c_2):b_R(b_2^{-1}c_1)\in B'\}=\{(b_1 c_1 , b_2 c_2):b_2^{-1}c_1\in C A\}
\end{equation}
%%%%%%%%
\notka{dpsi1}\begin{equation}\label{dpsi1}
D(\Psi_1):=\{(g_1,g_2) : b_L(g_2)\in B'\}=%\{(g_1,g_2) : s_B(g_2)\in BC\}=$$
%$$=\{(g_1,g_2) : s_B(b_L(g_2)c_R(g_1)^{-1})\in BC\}=
G_B\times b_L^{-1}(B')
\end{equation}
%%%%%%%%
\notka{fi1}\begin{equation}\label{fi1}
\Phi_1:D(\Phi_1) \ni(b_1 c_1, b_2 c_2)\mapsto 
(b_1 b_2 c_0, b_R(b_2 c_0)c_0^{-1} c_2)\,,\,\,c_0:=\tilde{c}_L(b_2^{-1}c_1)\end{equation}
\begin{equation}\label{psi1}\Psi_1: D(\Psi_1)\ni (b_1 c_1,b_2 c_2)\mapsto 
(b_1 b_R(b_2 c_1^{-1})^{-1} c_0, b_R(b_2 c_1^{-1}) c_1 c_2)\,,\,\,c_0:=\tilde{c}_L(s_B(b_2 c_1^{-1}))
\end{equation}
%%%%%%%%
\notka{dfi2}\begin{equation}\label{dfi2}D(\Phi_2):=\{(b_1 c_1 , b_2 c_2):b_1 b_2 c_2, b_2 c_2 \in C A \}\end{equation}
\notka{dpsi2}\begin{equation}\label{dpsi2}D(\Psi_2):=\{(b_1 c_1 ,b_2 c_2 ) : b_1 c_1 , b_2\in C A\}=
%\{(b_1 c_1 ,b_2 c_2 ) : b_1 c_1 , s_B(b_2 c_2)\in BC\}
b_R^{-1}(B')\times b_L^{-1}(B')\subset D(\Psi_1)
\end{equation}
\notka{fi2}\begin{equation} \label{fi2}
\Phi_2:D(\Phi_2) \ni(b_1 c_1, b_2 c_2)\mapsto 
(b_1 b_2 c_2, b_R(b_1 c_0) c_0^{-1} c_1)\,,\,\,c_0:=\tilde{c}_L(b_2 c_2) 
\end{equation}
$$\Psi_2: D(\Psi_2)\ni (b_1 c_1, b_2 c_2)\mapsto 
( c_0 \tilde{c_0}^{-1} b_2 c_2 , b_L^{-1} (c_0 \tilde{c_0}^{-1} b_2) b_1 c_1)\,,\,\,
c_0:=\tilde{c}_L(b_1 c_1)\,,\,
\tilde{c_0}:=\tilde{c}_L(b_2)$$

\begin{lem} Sets $D(\Phi_1), D(\Psi_1), D(\Phi_2), D(\Psi_2)$ are open.
$\Phi_1$ is a diffeomorphism of $D(\Phi_1)$ and $D(\Psi_1)$ and $\Psi_1=\Phi_1^{-1}$.
$\Phi_2$ is a diffeomorphism of $D(\Phi_2)$ and $D(\Psi_2)$ and $\Psi_2=\Phi_2^{-1}$.
\end{lem}
{\em Proof:} Short inspection of definitions shows that really these sets are open  and mappings are smooth. 
So it is sufficient to prove that $\Phi_1, \Psi_1$ and $\Phi_2, \Psi_2$ are pairs of mutually inverse mappings. 
The proof will be given for the first pair, for the second one it is similar.\\
1) The inclusion $\Phi_1(D(\Phi_1))\subset D(\Psi_1)$. 
Let $(b_1 c_1, b_2 c_2)\in D(\Phi_1)\,,\,b_2^{-1} c_1=:c_0 a_0$. 
$$\Phi_1(b_1 c_1, b_2 c_2)=(b_1 b_2 c_0, b_R(b_2 c_0)c_0^{-1} c_2)
=:(\tilde{b}_1\tilde{c}_1 , \tilde{b}_2\tilde{c}_2).$$
And we have: $s_B(\tilde{b}_2\tilde{c}_1^{-1})=s_B(b_R(b_2 c_0)c_0^{-1})=b_2 c_0=c_1a_0^{-1}\in C A$.\\
%%%%
2) The inclusion $\Psi_1(D(\Psi_1))\subset D(\Phi_1)$.
Let $(b_1 c_1, b_2 c_2)\in D(\Psi_1)\,,\,s_B(b_2  c_1^{-1})=b_R(b_2  c_1^{-1}) c_1=:c_0 a_0$. 
$$\Psi_1(b_1 c_1, b_2 c_2)=(b_1 b_R( b_2  c_1^{-1})^{-1} c_0, b_R( b_2  c_1^{-1})c_1 c_2)
=:(\tilde{b}_1\tilde{c}_1 , \tilde{b}_2\tilde{c}_2).$$
So $\tilde{b}_2^{-1} \tilde{c}_1=b_R( b_2  c_1^{-1})^{-1}c_0=c_1 a_0^{-1}\in C A$\\
3) Composition $\Psi_1 \Phi_1$:
$$\Psi_1\Phi_1(b_1 c_1, b_2 c_2)=\Psi_1(b_1 b_2 c_0, b_R(b_2 c_0) c_0^{-1} c_2)=$$
$$=\Psi_1(\tilde{b}_1\tilde{c}_1 , \tilde{b}_2\tilde{c}_2)=
(\tilde{b}_1 b_R(\tilde{b}_2\tilde{c}_1^{-1})^{-1}\tilde{c}_0, 
b_R(\tilde{b}_2\tilde{c}_1^{-1}) \tilde{c}_1 \tilde{c}_2)\,,\,\,{\rm where} \, 
s_B(\tilde{b}_2\tilde{c}_1^{-1})=:\tilde{c}_0\tilde{a}_0.$$
But $s_B(\tilde{b}_2\tilde{c}_1^{-1})=s_B(b_R(b_2 c_0) c_0^{-1})=b_2 c_0$ i.e. $\tilde{c}_0=c_1$, 
so we have:
$$\tilde{b}_1 b_R(\tilde{b}_2\tilde{c}_1^{-1})^{-1}\tilde{c}_0=
b_1 b_2 b_R(b_R(b_2 c_0)c_0^{-1})^{-1} c_1=b_1 c_1$$
$$b_R(\tilde{b}_2\tilde{c}_1^{-1}) \tilde{c}_1 \tilde{c}_2=b_2 c_0 c_0^{-1} c_2=b_2 c_2$$
4) Composition $\Phi_1 \Psi_1$:

$$\Phi_1 \Psi_1(b_1 c_1, b_2 c_2)=\Phi_1(b_1 b_R(b_2 c_1^{-1})^{-1} c_0, b_R(b_2 c_1^{-1}) c_1 c_2)=$$
$$=\Phi_1(\tilde{b}_1\tilde{c}_1 , \tilde{b}_2\tilde{c}_2)=
(\tilde{b}_1 \tilde{b}_2 \tilde{c}_0, b_R(\tilde{b}_2\tilde{c}_0)\tilde{c}_0^{-1}\tilde{c}_2),$$
where $s_B(b_2 c_1^{-1})=:c_0 a_0$ and $\tilde{b}_2^{-1} \tilde{c}_1=:\tilde{c}_0\tilde{a}_0$. 
But as shown in point 2) $\tilde{c}_0=c_1$, so
$$\tilde{b}_1 \tilde{b}_2 \tilde{c}_0=b_1 b_R(b_2 c_1^{-1})^{-1} b_R(b_2 c_1^{-1}) c_1=b_1 c_1$$
$$b_R(\tilde{b}_2\tilde{c}_0)\tilde{c}_0^{-1}\tilde{c}_2=
b_R(b_R(b_2 c_1^{-1})c_1) c_1^{-1} c_1 c_2=b_2 c_2$$
\dowl

By the lemma $\Phi_1^*:\sD(D(\Psi_1))\rightarrow \sD(D(\Phi_1))$  and 
$\Phi_2^*:\sD(D(\Psi_2))\rightarrow \sD(D(\Phi_2))$ are isomorphisms of vector spaces.

\begin{lem}\label{density-0}\notka{density-0}
There exists smooth, positive  function $k_1: D(\Phi_1)\rightarrow \R$ such that
$$\left(k_1 \Phi_1^*(f_1\mt f_2)\right)*(f_3\mt f_4)=\hat{\delta}_0(f_1) T^{-1}(f_3\mt f_2*f_4)$$
for $f_1,f_3,f_4\in \sD(G_B), f_2\in\sD(b_L^{-1}(B'))$.\\
There exists smooth, positive  function $k_2: D(\Phi_2)\rightarrow \R$ such that
$$\left(k_2 \Phi_2^*(f_1\mt f_2)\right)*(f_3\mt f_4)=\hat{\delta}_0(f_1) T^{-1}(f_2*f_3\mt f_4)$$
for $f_3\in \sD(G_B), f_2, f_4\in\sD(b_L^{-1}(B')), f_1\in\sD(b_R^{-1}(B'))$.
\end{lem}

Since in both cases the set of possible $f_3\mt f_4$ is linearly dense in $C_r^*(G_B)\mt C_r^*(G_B)$, the lemma means that
$$k_1 \Phi_1^*(f_1\mt f_2)=\hat{\delta}_0(f_1) T^{-1}(I\mt f_2)\,,\,f_1,\in \sD(G_B), f_2\in\sD(b_L^{-1}(B'))$$
$$k_2 \Phi_2^*(f_1\mt f_2)=\hat{\delta}_0(f_1) T^{-1}(f_2\mt I)\,,\,f_1\in\sD(b_R^{-1}(B')), f_2\in\sD(b_L^{-1}(B'))$$
and by density of subspaces $\sD(G_B),\sD(b_L^{-1}(B')),\sD(b_R^{-1}(B'))$ in $C_r^*(G_B)$ we obtain 
$$\Delta(C_r^*(G_B))(I\mt C_r^*(G_B)),\,\Delta(C_r^*(G_B))(C_r^*(G_B)\mt I)\subset C_r^*(G_B)\mt C_r^*(G_B).$$
{\em Proof of the lemma:} Let $f_1,f_2,f_3,f_4$ be as stated. Let us compute:
%%%%%%%
$$[k_1\Phi_1^*(f_1\mt f_2)]*(f_3\mt f_4)(b_1c_1,b_2c_2)=\int_{C\times C}d_l c\, d_l \tilde{c}\,k_1(b_1 c, b_2 \tilde{c})\times$$
$$\times  \, [\Phi_1^*(f_1\mt f_2)](b_1 c, b_2 \tilde{c})\, f_3(c_L^{-1}(b_1 c) b_1 c_1)\, f_4(c_L^{-1}(b_2 \tilde{c}) b_2 c_2)$$
Since $\Phi_1^*(f_1\mt f_2)\in \sD(D(\Phi_1))$ the integral over $c$ can be restricted to the set $C_{b_2^{-1}}:=\{c\in C: b_R(b_2^{-1}c)\in B'\}$.
Now, the map 
$\displaystyle C_{b_2^{-1}}\times C\ni (c,\tilde{c})\mapsto (c':=\tilde{c}_L(b_2^{-1}c), \tilde{c})\in C_{b_2}\times C$ 
is a diffeomorphism, as well as the map
$\displaystyle C_{b_2}\times C\ni (c',\tilde{c})\mapsto (c', c'':=c'^{-1}\tilde{c})\in C_{b_2}\times C.$
Their composition 
$$ C_{b_2^{-1}}\times C \ni \left(c,\tilde{c})\mapsto (c'=\tilde{c}_L(b_2^{-1}c), c''=\tilde{c}_L^{-1}(b_2^{-1}c)\tilde{c}\right)\in C_{b_2}\times C$$
has the inverse:
$$\Lambda_{b_2}: C_{b_2}\times C \ni (c',c'')\mapsto (c=\tilde{c}_L(b_2c'), \tilde{c}=c'c'')\in C_{b_2^{-1}}\times C$$
Using this change of variables we write the integral as:
\notka{calka1}\begin{equation}\label{calka1}
\int_{C_{b_2}\times C}d_lc'\, d_lc''\,|\det\Lambda_{b_2}'| \,k_1(b_1 \tilde{c}_L(b_2 c'), b_2 c'c'')\times$$
$$\times f_1(b_1 b_2 c')\,f_2(b_R(b_2c') c'')\, f_3(c_L^{-1}(b_1\tilde{c}_L(b_2 c')b_1 c_1)\, f_4 (c_L^{-1}(b_2 c' c'')b_2 c_2),
\end{equation}
moreover $|\det\Lambda_{b_2}'|$ is a smooth function of $b_2$.

On the other hand, using (\ref{T-form}) and (\ref{delta0-form}) we compute:
$$(\widehat{\delta_0}(f_1)T^{-1}(f_3\mt f_2* f_4))(b_1 c_1, b_2 c_2)=$$
$$=\int_C d_lc'\, j_C(c_L(b_2 c'))^{-1/2} f_1(b_1 b_2 c')
[T^{-1}(f_3\mt f_2* f_4)](c_L(b_1 b_2 c')^{-1} b_1 c_1, b_R(b_2 c')c^{-1} c_2)=$$
$$=\int_C d_lc'\, j_C(c_L(b_2 c'))^{-1/2} j_C(\tilde{c}_L(b_R(b_2 c')))^{-1/2}\, f_1(b_1 b_2 c')\times$$
$$\times f_3(c_L(b_1c_L(b_2 c')\tilde{c}_L(b_R(b_2 c')))^{-1} b_1 c_1)\, (f_2*f_4)(b_R(b_2 c')c'^{-1} c_2)$$
Using assumption about $f_2$ we can restrict domain of integration to the set $C_{b_2}$, 
and using formula for multiplication (\ref{mult-ga}) we expand the integral as:
$$=\int_{C_{b_2}\times C}d_lc' \,d_lc''\,
j_C(c_L(b_2 c'))^{-1/2} j_C(\tilde{c}_L(b_R(b_2 c')))^{-1/2}\, f_1(b_1 b_2 c') f_2(b_R(b_2 c') c'')\times$$ 
$$\times f_3(c_L(b_1c_L(b_2 c') \tilde{c}_L(b_R(b_2 c')))^{-1} b_1 c_1) f_4(c_L(b_2 c'c'')^{-1} b_2 c_2)$$
but $c_L(b_2 c') \tilde{c}_L(b_R(b_2 c')) =\tilde{c}_L(b_2 c')$ so we get:
$$=\int_{C_{b_2}\times C}d_lc' \,d_lc''\,
 j_C(\tilde{c}_L(b_2 c'))^{-1/2}\,
f_1(b_1 b_2 c') f_2(b_R(b_2 c')c'') \times$$
$$\times f_3(c_L(b_1 \tilde{c}_L(b_2 c'))^{-1} b_1 c_1) f_4(c_L(b_2 c' c'')^{-1} b_2 c_2)$$
Comparison of this expresion with (\ref{calka1}) gives the function $k_1$ that appears in the lemma.
%%%%%%%%%%%%%%%%%%%%%%%%%%%%%%%%%%%%%%%%%%%%%%%%%%%%%%%%%%%%%%%%%%%%%%%%%%%%%%%%%%%%%%%%%%%
%%%%%%%%%%%%%%%%%%%%%%%%%%%%%%%%%%%%%%%%%%%%%%%%%%%%%%%%%%%%%%%%%%%%%%%%%%%%%%%%%%%%%%%%%%

Now we pass to the second statement of the lemma \ref{density-0}. For $f_1,f_2,f_3,f_4$ as in the lemma let's compute:
$$[k_2\Phi_2^*(f_1\mt f_2)]*(f_3\mt f_4)(b_1c_1,b_2c_2)=\int_{C\times C}d_l c\, d_l \tilde{c}\,k_2(b_1 c, b_2 \tilde{c})\times$$
$$\times  [\Phi_2^*(f_1\mt f_2)](b_1 c, b_2 \tilde{c})\, f_3(c_L^{-1}(b_1 c) b_1 c_1)\, f_4(c_L^{-1}(b_2 \tilde{c}) b_2 c_2)$$
Since $\Phi_2^*(f_1\mt f_2)\in \sD(D(\Phi_2))$ the integral over $\tilde{c}$ can be restricted to the set $C_{b_2}$, 
and using formula (\ref{fi2}) we obtain:
$$\int_{C\times C_{b_2}}d_l c\, d_l \tilde{c}\, k_2(b_1 c, b_2 \tilde{c})\,f_1(b_1 b_2 \tilde{c})\times$$
$$\times \, f_2(b_R(b_1\tilde{c}_L(b_2\tilde{c}))\tilde{c}_L^{-1}(b_2\tilde{c}) c)\,
f_3(c_L^{-1}(b_1 c) b_1 c_1)\, f_4(c_L^{-1}(b_2 \tilde{c}) b_2 c_2)$$
The mapping 
$$C\times C_{b_2}\ni (c,\tilde{c})\mapsto (c':=\tilde{c}_L^{-1}(b_2 \tilde{c}) c, c'':=\tilde{c})\in C\times C_{b_2}$$
is a diffeomorphism with the inverse:
$$\tilde{\Lambda}_{b_2}:C\times C_{b_2}\ni (c',c'')\mapsto (c:=\tilde{c}_L(b_2 c'') c', \tilde{c}:=c'')\in C\times C_{b_2}$$
Therefore our integral is equal to:
\begin{equation}\label{calka2}
\int_{C\times C_{b_2}}d_l c'\, d_l c''\,|\det\tilde{\Lambda}_{b_2}'|\, k_2(b_1 \tilde{c}_L(b_2c'')c' , b_2 c'')\, f_1(b_1 b_2 c'')\times$$
$$\times \, f_2(b_R(b_1\tilde{c}_L(b_2 c'')) c')\,
f_3(c_L^{-1}(b_1 \tilde{c}_L(b_2c'') c') b_1 c_1)\, f_4(c_L^{-1}(b_2 c'') b_2 c_2),
\end{equation}
and again $|\det\tilde{\Lambda}_{b_2}'|$ is a smooth function of $b_2$.
%%%%%%%%%%

On the other hand using (\ref{T-form}) and (\ref{delta0-form}) we compute: 
$$(\widehat{\delta_0}(f_1)T^{-1}(f_2* f_3\mt f_4))(b_1 c_1, b_2 c_2)=$$
$$=\int_C d_lc''\, j_C(c_L(b_2 c''))^{-1/2} j_C(\tilde{c}_L(b_R(b_2 c'')))^{-1/2}\, f_1(b_1 b_2 c'')\times$$
$$\times (f_2* f_3)(c_L(b_1c_L(b_2 c')\tilde{c}_L(b_R(b_2 c')))^{-1} b_1 c_1)\, f_4(b_R(b_2 c')c'^{-1} c_2).$$
Again  using the assumption about the support of $f_4$ and expanding the product $f_2 *f_3$ 
we can write this integral as:
$$\int_{C_{b_2}\times C} d_l c''\,d_l c'\,j_C(\tilde{c}_L(b_2 c''))^{-\frac12} 
f_1(b_1 b_2 c'') f_2(b_R(b_1\tilde{c}_L(b_2 c''))c')\times$$
$$\times f_3(c_L(b_1 \tilde{c}_L(b_2 c'')c')^{-1}b_1 c_1)\, f_4(c_L(b_2 c'')^{-1} b_2 c_2)$$
Again comparison of this expresion with (\ref{calka2}) gives the function $k_2$ that appears in the lemma.

\dowl

Finally, from the previous lemma and lemma \ref{ind-lim} to prove density conditions it is sufficient to prove the following:
\begin{lem}\label{gestosc}\notka{gestosc}
The linear spaces $\sA(D(\Phi_1))$ and $\sA(D(\Phi_2))$ are dense in
$C^*_r(G_B)\mt C^*_r(G_B)$.
\end{lem}
{\em Proof:}  It is sufficient to show that a function $F\in \sD(G_B\times
G_B)$ can be approximated in norm $||\cdot||_0$ by functions supported
in $D(\Phi_1)$ and a function  $F\in\sD(G_B\times b_R^{-1}(B'))$ (by \ref{zal} (3)) by functions supported
in $D(\Phi_2)$.

Let us begin with $D(\Phi_1)$.
For an open  $V\subset B$ such that $B\setminus B'\subset V$  let $\tilde{\chi}$ be a smooth function on $A$ 
satisfying conditions: $\tilde{\chi}=1$ on some neighbourhood of $B\setminus B'$; $0\leq\tilde{\chi}\leq 1$
and $supp ( \tilde{\chi})\subset V$. Define:
$$\chi(b_1 c_1, b_2 c_2):=\tilde{\chi}(b_R(b_2^{-1} c_1))$$
Let $F\in\sD(G_B\times G_B)\,,\,supp(F)=:K_F$, then $F=(F-F \chi)+F \chi$ and 
$(F-F \chi)\in\sD(D(\Phi_1))$ so we have to prove that, by choosing $V$, the $C^*$-norm of $F \chi$ 
can be made as small as we wish. It is sufficient to prove that for norms $||F \chi||_l$ and 
$||F \chi||_r=||(F \chi)^*||_l$.
For the left norm we need to estimate  (compare (\ref{l-norm})) the integral
$$I_l:=\int d_lc_1 d_lc_2 |F\chi(b_1 c_1, b_2 c_2)|\,\,\,{\rm for}\,\,\,(b_1, b_2)\in (b_L\times b_L)(K_F).$$
First we estimate this integral from above by $ \sup|F| \int_{C(b_1, b_2)} d_lc_1 d_lc_2$, where 
$$C(b_1,b_2):=\{(c_1,c_2): (b_1 c_1, b_2 c_2)\in K_F\,,\, b_R(b_2^{-1} c_1)\in supp (\tilde{\chi})\};$$
Now we have the chain of inclusions:
$$C(b_1,b_2)=\{(c_1,c_2): (b_1 c_1,b_2 c_2)\in K_F\}\cap \{(c_1,c_2): b_R(b_2^{-1} c_1)\in supp (\tilde{\chi})\}\subset$$
$$\subset (c_R\times c_R)(K_F)\cap(\{c\in C: b_R(b_2^{-1}c)\in supp (\tilde{\chi})\}\times C)\subset$$
$$\subset (K_C\times K_C)\cap (\{c\in C: b_R(b_2^{-1}c)\in supp (\tilde{\chi})\}\times C)=$$
$$=(K_C\cap \{c\in C: b_R(b_2^{-1}c)\in supp (\tilde{\chi})\})\times K_C\subset  Z(b_2^{-1},e,K_C; V)\times K_C$$
where $K_C\subset C$ is a compact set such that $(c_R\times c_R)(K_F)\subset K_C\times K_C$; note that $K_C$ depends only on $F$.
In this way we get the estimate for $(b_1, b_2)\in (b_L\times b_L)(K_F)$:
$$I_l\leq \sup|F| \left(\int_{K_C} d_lc\right) \mu(b_2^{-1},e, K_C; V)$$
and, if $K_1$ is a compact subset of $B$ such that $(b_L\times b_L)(K_F)\subset K_1\times K_1$, we obtain
\notka{lewanorma1}\begin{equation}\label{lewanorma1} ||F\chi||_l\leq  \sup|F| \left(\int_{K_C} d_lc\right) \mu(K_1^{-1},\{e\}, K_C; V)
\end{equation}
%%%%%%%%%%%%%
Now the right norm of $F\chi$  i.e. the left norm of $(F\chi)^*$. The integral to estimate is:
$$I_r:=\int d_lc_1 d_lc_2 |(F\chi)^*(b_1 c_1, b_2 c_2)|$$
As above, first we estimate this integral by $\displaystyle \sup|F|\int_{\tilde{C}(b_1,b_2)} d_lc_1 d_lc_2$, 
where 
$$\tilde{C}(b_1,b_2):=\{(c_1,c_2): (b_1 c_1, b_2 c_2)\in supp(F\chi)^*\}.$$
 Now,
$$(b_1 c_1, b_2 c_2)\in supp(F\chi)^*\iff(s_B(b_1 c_1), s_B(b_2 c_2))\in supp(F\chi)\iff$$
$$\iff\left[(s_B(b_1 c_1), s_B(b_2 c_2))\in K_F \wedge b_R(b_R(b_2 c_2)^{-1}c_1^{-1})\in supp(\tilde{\chi})\right]$$
Let us denote $\tilde{K}_F:=(s_B\times s_B)(K_F)$, so  we have to estimate the integral 
$\displaystyle\int_{\tilde{C}(b_1,b_2)} d_lc_1 d_lc_2$ for $(b_1,b_2)\in(b_L\times b_L)(\tilde{K}_F)$. Again there is an inclusion:
$$\tilde{C}(b_1,b_2)\subset (\tilde{K}_C\times \tilde{K}_C)\cap \{(c_1, c_2): b_R(b_R(b_2 c_2)^{-1}c_1^{-1})\in supp \tilde{\chi} \},$$
where $\tilde{K}_C$ is a compact set such that $(c_R\times c_R)(\tilde{K}_F)\subset \tilde{K}_C\times \tilde{K}_C$.\\
For fixed $c_2\in\tilde{K}_C$ define  $\tilde{b}:=b_R(b_2 c_2)$ and consider the set 
$$\{c\in\tilde{K}_C: b_R(\tilde{b}^{-1}c^{-1})\in supp(\tilde{\chi})\}=
\{c\in\tilde{K}_C^{-1}: b_R(\tilde{b}^{-1}c)\in supp(\tilde{\chi})\}^{-1}\subset Z(\tilde{b},e, \tilde{K}_C^{-1},V)^{-1}$$
 Since all the sets $Z(\tilde{b},e,\tilde{K}_C^{-1};V)$ are subsets of the fixed compact set $\tilde{K}_C^{-1}$ there exists a constant 
$M(\tilde{K}_C)$ depending 
only on $\tilde{K}_C$, such that 
$$\int_{Z(\tilde{b},e,\tilde{K}_C^{-1};V)^{-1}} d_lc\leq M(\tilde{K}_C) \int_{Z(\tilde{b},e,\tilde{K}_C^{-1};V)} d_lc=
M(\tilde{K}_C)\mu(\tilde{b},e,\tilde{K}_C^{-1};V)$$
But $\tilde{b}=b_R(b_2 c_2)$ belongs to a fixed compact set $K_2=b_R(\tilde{K}_B\tilde{K}_C)$, where 
$\tilde{K}_B$ is a compact set such that $(b_L\times b_L)\tilde{K}_F\subset \tilde{K}_B\times \tilde{K}_B$
Putting all together we get the estimate:
\notka{prawanorma1}\begin{equation}\label{prawanorma1}
 ||F\chi||_r\leq \sup |F| M(\tilde{K}_C)\mu(K_2,\{e\},\tilde{K}_C^{-1};V)\left(\int_{\tilde{K}_C} d_lc \right)
\end{equation}
This inequality together with (\ref{lewanorma1}) and the condition (4) of assumptions (\ref{zal}) gives the density of $\sA(D(\Phi_1))$. 

Now we pass to $D(\Phi_2)$. Let $\tilde{\chi}$ be as above and define 
$\chi(b_1 c_1, b_2 c_2):=\tilde{\chi}(b_R(b_1 b_2 c_2))$. It is sufficient to approximate functions 
$F\in\sD(G_B\times b_R^{-1}(B'))$ and to this end we need an estimate for $||F\chi||_l$
and $||(F\chi)^*||_l$. 
Let $K_F:=supp(F)$ and $(b_1, b_2)\in (b_L\times
b_L)(K_F)$. As above we have inequality
$$\int d_lc_1 d_lc_2 |F\chi| \leq \sup|F|\int_{C(b_1,b_2)} d_lc_1 d_lc_2,\,{\rm where}$$ 
$$C(b_1,b_2):=\{(c_1,c_2): (b_1 c_1, b_2 c_2)\in K_F\,,\,b_R(b_1 b_2 c_2)\in supp(\tilde{\chi})\}\subset$$
$$\subset (K_C\times K_C)\cap (C\times \{c\in C: b_R(b_1 b_2 c)\in supp(\tilde{\chi})\})=$$
$$=K_C\times(K_C\cap \{c\in C: b_R(b_1 b_2 c)\in supp(\tilde{\chi})\})=K_C\times Z(b_1 b_2,e, K_C,V),$$
where $K_C\subset C$ is a compact such that $(c_R\times c_R)(K_F)\subset K_C\times K_C$;
In this way we obtain the estimate:
\notka{lewanorma2}\begin{equation}\label{lewanorma2}
 ||F\chi||_l\leq \sup |F| \left(\int_{K_C} d_lc\right)\mu(K_1,\{e\},K_C;V),\, where \end{equation}
$K_1=K_B K_B$  for a compact $K_B\subset B$ with $(b_L\times b_L)(K_F)\subset K_B\times K_B$.

Finally we pass to  $(F\chi)^*$. As above, we consider the set
$$\tilde{C}(b_1,b_2):=\{(c_1,c_2): (s_B(b_1 c_1), s_B(b_2 c_2))\in supp(F \chi)\}$$
Let $\tilde{K}_F:=(s_B\times s_B)(K_F)$ and $\tilde{K}_C\subset C$ be a compact set such that
$(c_R\times c_R)(\tilde{K}_F)\subset \tilde{K}_C\times \tilde{K}_C$. It is sufficient to consider 
$(b_1,b_2)\in (b_L\times b_L)(\tilde{K}_F)$.
$$\tilde{C}(b_1,b_2)\subset(\tilde{K}_C\times \tilde{K}_C)\cap\{(c_1,c_2): b_R(b_R(b_1 c_1)c_L(b_2 c_2)^{-1})b_2\in V\}$$
Let us fix $c_1\in \tilde{K}_C$, denote $b_R(b_1 c_1)=:\tilde{b}$ and consider the set
$$S(b_1,b_2, c_1):=\tilde{K}_C\cap\{c\in C: b_R(\tilde{b}c_L(b_2 c)^{-1})b_2\in V\}.$$
The map $\varphi_{b_2}:c\mapsto c_L(b_2 c)$ is a diffeomorphism so
$$\{c\in C: b_R(\tilde{b}c_L(b_2 c)^{-1})b_2\in V\}=\varphi_{b_2}^{-1}(\{c\in C: b_R(\tilde{b}c^{-1})b_2\in V\})$$
and
$$\tilde{K}_C\cap\{c\in C: b_R(\tilde{b}c_L(b_2 c)^{-1})b_2\in V\}=
\varphi_{b_2}^{-1}\left(\varphi_{b_2}(\tilde{K}_B)\cap \{c\in C: b_R(\tilde{b}c_L(b_2 c)^{-1})b_2\in V\}\right)$$
The sets $\varphi_{b_2}(\tilde{K}_C)=c_L(b_2 \tilde{K}_C)$ are contained in a fixed compact set $H_C\subset C$ therefore
$$\varphi_{b_2}^{-1}\left(\varphi_{b_2}(\tilde{K}_C)\cap \{c\in C: b_R(\tilde{b}c_L(b_2 c)^{-1})b_2\in V\}\right)\subset
\varphi_{b_2}^{-1}\left(\{c\in H_C: b_R(\tilde{b}c^{-1})b_2\in V\}\right)=$$
$$=\varphi_{b_2}^{-1}\left(Z(\tilde{b},b_2,H_C^{-1};V)^{-1}\right)$$
Since $\varphi_{b_2}^{-1}=\varphi_{b_2^{-1}}$ and $b_2$ is contained in a fixed compact set depending only on
$\tilde{K}_F$ there exists a constant $W(\tilde{K}_F)$ such that
$$\int_{\varphi_{b_2}^{-1}(K)} d_lc\leq W(\tilde{K}_F)\int_K d_lc$$ in this way, for fixed $c_1$ we obtain that
$$\int_{S(b_1,b_2,c_1)} d_lc_2\leq W(\tilde{K}_F) M(H_C)\mu(\tilde{b},b_2, H_C^{-1};V),$$
where $M(H_B)$ is as before (\ref{prawanorma1}). Finally, since $\tilde{b}=b_R(b_1 c_1)$ is contained in a 
compact set $\tilde{K}_B$ depending only on $K_F$
and $b_2$ is contained in a fixed compact $\tilde{K}_1\subset B'$ we obtain the estimate
$$||(F\chi)^*||_l\leq  \sup|F|\, W(\tilde{K}_F) M(H_B) \mu(\tilde{K}_B, \tilde{K}_1, H_B^{-1};V)\left(\int_{\tilde{K}_B} d_lc\right)$$
This estimate together with (\ref{lewanorma2}) and condition (4) of assumptions (\ref{zal}) proves the density of $\sA(D(\Phi_2))$.\\
\dowl

This completes the proof of the statement b) of prop \ref{Delta}.\\\dow
% \begin{center}{\large \bf Fakty do warunkow gestosci:}\end{center}
%%%%%%%%%%%%%%%%%%%%%%%%%%%%%%%%%%%%%%%%%%%%%%%%%%%%%%%%%%%%%%%%%%%%%%%%%%%%%%%%%%%%%%%%%%%%%%%%%
%%%%%%%%%%%%%%%%%%%%%%%%%%%%%%%%%%%%%%%%%%%%%%%%%%%%%%%%%%%%%%%%%%%%%%%%%%%%%%%%%%%%%%%%%%%%%%%%%%%%
\section{'ax+b' again}

Now  we apply results of the previous section to the 'ax+b' group. In this section 
$G$ is the 'ax+b' group; $$B:=\{(b,1)\,,\,b\in\R\}\,,\,\,C:=\{(c-1,c)\,,\,c\in\R_*\}\,,\,\,A:=\{(0,a)\,,\,a\in\R_*\}$$
We have $A\cap C =B\cap C=\{e\}$ and $G=B C$.  

% We apologize the reader for changing roles of subgroup and present 
% the short 'dictionary' from section \ref{stwist} to 'ax+b' notation.

% $$\renewcommand{\arraystretch}{1.3}  
% \begin{array}{lrccccccccc}
% {\rm section}  \,\ref{stwist} &\left|\right. & A & B & C& A'=A\cap B C  & a_L, a_R & b_L, b_R & \tilde{b}_L, \tilde{b}_R &
% \delta_0=m_B^T & \tilde{m}_B^T\\\hline
% {\rm ax+b}  &\left| \right.& B & C & A & B'=B\cap C A &  b_L, b_R & c_L, c_R & \tilde{c}_L, \tilde{c}_R & \delta_0=m_C^T & \tilde m_C^T
% \end{array}
% \renewcommand{\arraystretch}{0.9}
% $$

So our main object here is a groupoid $G_B$ related to the double Lie group $(G;B,C)$ together with coassociative morphism 
$\delta_0=m_C^T: G_B\rel G_B\times G_B$. We identify $G$ with $B\times C=\R\times \R_*$ by:
$$\R\times\R_*\ni(z,c)\mapsto (z,1)(c-1,c)=(z+c-1,c)\in G$$
In this presentation the relevant objects are given by:
$$b_L(z,c)=(z,1)\,,\,b_R(z,c)=(\frac{z}{c},1)\,,\,c_L(z,c)=c_R(z,c)=(0,c)\,,\,\tilde{c}_L(z,c)=(0,z+c)$$
$$m_B=\{(z,c_1 c_2; z, c_1,\frac{z}{c_1},c_2)\}\,\,,\,\,\delta_0=\{(z_1,c,z_2,c;z_1+z_2,c)\}\,\,,\,\,
B'=\{(z,1) : z+1\neq 0\}$$
$$T:=\{(z_1,\frac{1}{1+z_2},z_2,1): z_2+1\neq 0\}$$
$$\tilde{m}_C:=\{(z_1+z_2+z_1 z_2,c; z_1,\frac{c+z_2}{1+z_2},z_2,c) : (1+z_1)(1+z_2)(z_2+c)(z_1+z_2+z_1 z_2+c)\neq 0\}$$
%%%%%%%%
%
We will also need explicit formulae for operations in $\sA(G_B)$.

{\em The choice of $\om_0$.}  Let  $\lambda_0(z,1)(\partial_c):=1$. This choice leads to the left invariant half 
density $\lambda_0$,  the corresponding right invariant half density $\rho_0$ and $\om_0:=\lambda_0\mt\rho_0$. 
The formulae are:  
\begin{equation}
\lambda_0(z,c)(c  \partial_c)=\rho_0(z,c)(z \partial_z+c \partial_c)=1
\end{equation}
The  multiplication and $*$-operation  in $\sA(G_B)$ are given by: 
\begin{equation} 
(f_1*f_2)(z,c):=\int_{\R_*} \frac{{\rm d c_1} }{|c_1|} f_1(z,c_1) f_2(z c_1^{-1},c c_1^{-1})\,\,,\,
f^*(z,c):=\overline{f(z c^{-1}, c^{-1})}.
\end{equation}
%%%%%%%%%%%%%
{\em The choice of $\psi_0$.}  We choose $\nu_0(z):=| d z|^{1/2}$ and obtain 
$\psi_0(z,c)(\partial_z,\partial_c)=\frac{1}{|c|}$, and the formula for the scalar product in $\sD(G_B)$ reads:
$$(f_1,f_2)=\int_\Gamma \frac{{\rm d z}\,  {\rm d c}}{c^2} \overline{f_1(z,c)} f_2(z,c)$$
and the formula for $\pi_{id}$:
$$(\pi_{id}(f)\psi)(z,c)=\int_{\R_*} \frac{{\rm d c_1}}{|c_1|} f(z,c_1) \psi(z c_1^{-1},c c_1^{-1})\,,\,f\in\sA(G_B)\,,\,\,
\psi\in \sD(G_B)$$

Now we are going to verify assumptions \ref{zal}. The first and the second one are obvious, so let's pass to the third one.
We have to verify that $\sA(b_L^{-1}(B'))$ is dense in $C^*_r(G_B)$. Take  $f,\psi\in\sD(G_B)$ and compute:
$$||\pi_{id}(f)\psi||^2=\int \frac{{\rm dy}\,{\rm dc}}{c^2} |\pi_{id}(f)\psi(y,c)|^2=
\int \frac{{\rm dy}\,{\rm dc}}{c^2}\left|\int \frac{{\rm db}}{|b|}  f(y,b)\psi(yb^{-1},c b^{-1})\right|^2$$
Now we use Schwartz  inequality for functions: $b\mapsto f(y,b)$ and $b\mapsto\psi(yb^{-1},c b^{-1})$: 
$$\left|\int\frac{{\rm db}}{|b|}  f(y,b)\psi(yb^{-1},c b^{-1})\right|^2\leq 
\left(\int\frac{{\rm db}}{|b|}|f(y,b)|^2\right)\, \left(\int\frac{{\rm db}}{|b|}  |\psi(yb^{-1},c b^{-1})|^2\right)$$
and get the estimate:
$$||\pi_{id}(f)\psi||^2\leq \int {\rm dy} \left(\int\frac{{\rm db}}{|b|}  |f(y,b)|^2\right)
\int \frac{{\rm dc} \,  {\rm db} }{|b| c^2} |\psi(yb^{-1},c b^{-1})|^2$$
Note that for $y\neq 0$ applying the change of variables $(b,c)\mapsto (yb^{-1},c b^{-1})$ to the integral on the right 
we obtain:
$$\int\frac{{\rm dc} \,  {\rm db} }{|b| c^2}  |\psi(yb^{-1},c b^{-1})|^2=
\frac{1}{|y|}\int \frac{{\rm dc} \,  {\rm db} }{c^2}  |\psi(b,c)|^2$$
Therefore if $supp \,\psi\subset\{(z,c):z\neq0\}$ we get:
\notka{norm-est}\begin{equation}\label{norm-est}
||\pi_{id}(f)\psi||^2\leq ||\psi||^2 \left(\int \frac{{\rm dy}\, {\rm db}}{|y b|} |f(y,b)|^2\right)=:||\psi||^2||f||_2^2
\end{equation}
Since the set of such $\psi$'s  is dense the estimate is valid for any  $\psi$ and we get $||f||\leq ||f||_2$.

For  $\epsilon > 0$, let $\chi_{\epsilon}: \R\rightarrow [0,1] $ be a smooth function that is $1$ on some neighbourhood of $0$ 
and $0$ on the set $\{x\in \R :|x|\geq \epsilon\}$. Let $f\in \sD(G_B)$ and for $z_0\neq 0$ let 
$f_{\epsilon}(z,c):=\chi_{\epsilon}(z-z_0)f(z,c)$. Now the estimate given above implies 
$f_{\epsilon}\rightarrow 0$ in $C^*_r(G_B)$ as $\epsilon\rightarrow 0$. But $f=(f-f_\epsilon)+f_\epsilon$ and 
$f-f_\epsilon\in \sD(b_L^{-1}(\R\setminus\{z_0\}))$. This proves that statement (3) of (\ref{zal}) is true in our situation.

Now we are going to verify the fourth condition in assumptions (\ref{zal}) which in our situation takes form:\\
Let $K_C\subset C$ be compact,   $V\subset B$ open and $(z_1,z_2)\in B\times B'$.  Let 
$$Z(z_1,z_2,K_C;V):=K_C\cap\{c\in C: b_R(z_1, c)(z_2,1) \in V\},$$
and the function $\mu(z_1,z_2,K_C;V)$ be given by:
$$B\times B' \ni(z_1,z_2)\mapsto \mu(z_1,z_2,K_C;V):=\int_{Z(z_1,z_2,K_C;V)} d_l c.$$
Let  $K_1\subset B$ and $K_2\subset B'$ be compact and  $\mu(K_1,K_2,K_C;V):=\sup\{\mu(z_1,z_2,K_C;V): z_1\in K_1\,,\,z_2\in K_2\}$
Then \begin{equation}
\forall\,\epsilon>0\,\exists\, V-{\rm a\, neighbourhood\,of\,}B\setminus B'{\rm\, in\, B}\,:\, \mu(K_1,K_2,K_C;V)\leq\epsilon
\end{equation}

It is sufficient to check this condition for $K_C=K_m:=\{ c: m\leq |c|\leq \frac1m\}\,,\,m<1$. 
Using formula for $b_R$ we get $b_R(z_1,c)(z_2,1)=(\frac{z_1}{c}+z_2,1)$ and 
$$Z(z_1,z_2.K_m;V)=\{c\in\R_*: m\leq |c|\leq \frac1m\,,\,\frac{z_1}{c}+z_2\in V\}\,,\,1+z_2\neq 0.$$
The left invariant measure on $C$ is $\frac{dc}{|c|}$  and we obtain:
$$\mu(z_1,z_2,K_m;V):=\int_{Z(z_1,z_2,K_m;V)} \frac{dc}{|c|}\leq \frac1m \int_{Z(z_1,z_2,K_m;V)} dc$$
We will look for $V=V_\delta:=\{z\in \R: |z+1|<\delta\}$. In this situation $Z(z_1,z_2,K_m;V_\delta)$ is given by $c$ satisfying 
inequalities:
$$m\leq |c|\leq \frac1m\,,\,\,\left|\frac{z_1}{c}+z_2+1\right|<\delta\,\,\,,\,(1+z_2\neq 0)$$
The second inequality is for $\delta<|1+z_2|$ equivalent to:
$$\frac{-sgn(c) sgn(1+z_2) z_1}{\delta+|1+z_2|}<|c|<\frac{-sgn(c) sgn(1+z_2) z_1}{|1+z_2|-\delta}$$
and we get for $sgn(1+z_2) z_1\geq 0$:
$$c<0\,,\,m\leq |c|\leq \frac1m\,,\,\frac{|z_1|}{\delta+|1+z_2|}<|c|<\frac{|z_1|}{|1+z_2|-\delta}$$
and for $sgn(1+z_2) z_1<0 $:
$$c>0\,,\,m\leq |c|\leq \frac1m\,,\,\frac{|z_1|}{\delta+|1+z_2|}<|c|<\frac{|z_1|}{|1+z_2|-\delta}$$
Therefore the integral  $\int_{Z(z_1,z_2,K_m;V)} dc$ is majorized by:
$$\frac{|z_1|}{|1+z_2|-\delta}-\frac{|z_1|}{\delta+|1+z_2|}=\frac{2\delta |z_1|}{|1+z_2|^2-\delta^2}$$
and there is an estimate (for $\delta<|1+z_2|$):
$$\mu(z_1,z_2,K_m;V_\delta)\leq\delta  \frac{2 |z_1|}{m (|1+z_2|^2-\delta^2)}$$

It is sufficient to consider $K_1=K_M:=\{z\in\R: |z|\leq M\}$ and $K_2=\tilde{K}_M:=\{z\in \R: \frac1M\leq |z+1|\leq M\}$ for $M>1$.
For $(z_1,z_2)\in K_M\times \tilde{K}_M$ and $\delta\leq\frac1{2M}$ we have:
$$\mu(z_1,z_2,K_m;V_\delta)\leq\delta  \frac{2M}{m (|1+z_2|^2-\delta^2)}\leq \delta\frac{2M}{m(1/M^2-1/(4 M^2))}=\delta\frac{8 M^3}{3m}$$
So there is an estimate: $$\mu(K_M,\tilde{K}_M,K_m;V_\delta)\leq \delta  \frac{8 M^3}{3m}$$
and the fourth statement of assumptions (\ref{zal}) is fulfilled. Therefore we have
\begin{prop} Let $G$ be the 'ax+b' group and $A:=\{(0,a): a\in \R_*\}$, $B:=\{(b,1): b\in \R\}$,  $C:=\{(c-1,c): c\in \R_*\}$.
The reduced $C^*$ algebra of a differential groupoid $\Gamma=AC\cap CA$ over $A$ is 
isomorphic to $C_r^*(G_B)$ -- the reduced $C^*$ algebra of a differential groupoid related to a double Lie group $(G;B,C)$. 
It is equipped with a comultiplication $\Delta=\widehat{\sT}\Delta_0\widehat{\sT}^{-1}$ satysfying the density conditions, where 
$\widehat{\sT}$ is a unitary multiplier and satisfies: 
$$(\widehat{\sT}\mt I)(\Delta_0\mt id)\widehat{\sT}=(I \mt \widehat{\sT})(id\mt \Delta_0)\widehat{\sT}$$
\end{prop}

%%%%%%%%%%%%%%%%%%%%%%%%%%%%%%%%%%%%%%%%%%%%%%%%%%%%%%%%%%%%%%%%%%%%%%%%%%%%%%%%%%%%%%
%%%%%%%%%%%%%%%%%%%%%%%%%%%%%%%%%%%%%%%%%%%%%%%%%%%%%%%%%%%%%%%%%%%%%%%%%%%%%%%%%%%%%%
\section{The $C^*$-algebra and generators}
In this section we denote the groupoid $G_B$ by $\Gamma$ and identify it with the right transformation groupoid $\R\times\R_*$ with the action
$(z,c)\mapsto z/c$. The set of units will be denoted by $E$ and the set $\{z: z+1\neq 0\}$ by $E'$.
For $t\in\R$ let  $B_t:=\{(z,e^{t}): z\in \R\}\subset \Gamma$. Then $B_t$ is a one parameter group of bisections of $\Gamma$, 
therefore it defines one parameter group $\hat{B}_t$ of unitaries on $L^2(\Gamma)$ which are multipliers 
of $C^*_r(\Gamma)$. Easy computations show that $B_t\psi_0=\psi_0$ (where $\psi_0$ was defined in section 2)  so 
$\hat{B}_t(f\psi_0)=:(\hat{B}_t f)\psi_0$ and 
$$(\hat{B}_t f) (z,c)=f(e^{-t}z, e^{-t}c)\,,\,f\in\sD(\Gamma).$$

The set  $J:=\{(z,-1): z\in \R\}$ is a bisection, $J^2=E$,   so it defines unitary, selfadjoint operator $\hat{J}$
which is a multiplier of $C^*_r(\Gamma)$. Again $J\psi_0=\psi_0$ so 
$$(\hat{J}f)(z,c):= f (-z,-c)\,,\, f\in\sD(\Gamma)$$
Let $Y:\R\ni z \mapsto z \in \R$ and let the  operator  $\hat{Y}$ be defined by: 
$$(\hat{Y}f)(z,c):=Y(z) f(z,c)=z f(z,c) \,,\,f\in\sD(\Gamma)$$
Then the closure of $\hat{Y}$, which will be denoted by the same letter, is an unbounded, 
selfadjoint operator affiliated with  $C^*_r(\Gamma)$. 
On $\sD(\Gamma)$ these operators satisfy commutation relations: 
$$\hat{J}\hat{B}_t=\hat{B}_t\hat{J}\,\,,\,\,
\hat{J}\hat{Y}+\hat{Y}\hat{J}=0\,,\,\,e^{t}\hat{B}_t\hat{Y}=\hat{Y}\hat{B}_t$$
Note that the last relation gives precise meaning to the relation $[\hat{X},\hat{Y}]=i \hat{Y}$, 
for $\hat{X}$- generator of $\hat{B}_t$.
Additionaly if $\hat{Y}=sgn(\hat{Y})|\hat{Y}|$ is a polar decomposition then .
$$|\hat{Y}|\hat{J}=\hat{J}|\hat{Y}|\,\,,\,\, e^{t}\hat{B}_t |\hat{Y}|=|\hat{Y}| \hat{B}_t$$

\begin{lem} \label{op-X}\notka{op-X} The group $\hat{B}_t$ is strictly continuous. 
Let $\hat{X}$ denote the generator of $\hat{B}_t$, then $\hat{X}$ is 
affiliated with $C^*_r(\Gamma)$.  On $\sD(\Gamma)$ 
$\hat{X} f=i (z \partial_z+c\partial_c)f$ and $\sD(\Gamma)$ is a core for $\hat{X}$.
\end{lem}
{\em Proof:} Since $\hat{B}_t$ is a group of unitary multipliers, strict continuity is equivalent to continuity
of the mapping $\R\ni t\mapsto B_t f\in\sA(\Gamma)$ in $t=0$ for any $f\in\sA(\Gamma)$. 
The mapping $\Phi:\R\times\Gamma\ni (t,\gamma)\mapsto B_t(\gamma)\in \Gamma$ is continuous. Therefore if $|t|\leq \delta$ 
support of $B_tf$ is contained in a compact set $\Phi([-\delta,\delta]\times supp\,f)$ and $B_tf$ converges to $f$ uniformly 
as $t$ goes to $0$. So, by the lemma \ref{ind-lim} also in $C^*_r(\Gamma)$. By the general result \cite{Wor1}
$\hat{X}$ is affiliated to  $C^*_r(\Gamma)$. To prove the second claim it is enough to show that $\frac1t(B_tf-f)$ converges to
$s (z \partial_z+c\partial_c)f$ uniformly, but since $f$ is smooth this is clear. Finally, since $\sD(\Gamma)$ is 
$B_t$-invariant, it is a core for $\hat{X}$.\dowl

\begin{prop} $C^*_r(\Gamma)$ is generated by $\hat{X},\hat{Y},\hat{J}$.
\end{prop}
{\em Proof:} The proof is based on identification of $C^*_r(\Gamma)$ with crossed product $C_0(\R)\times_{\alpha}\R_*$. 
This is known, but we need explicit form of the isomorphism so we present relevant formulae. 
Let a Lie group $H$ acts from the right on a manifold $X$. Let $\alpha$ be the  corresponding (left) action 
on $C_0(X)$ ( i.e. $(\alpha_h f)(x):=f(x h)$). Since on a Lie group the Haar measure is given by a left invariant density, 
we choose such a real half density $\lambda$. Then we have $\int_H f(h) d_l(h)=\int_H f \lambda^2\,,\,f\in\sD(H)$. 
The modular function is given by $\delta(h)=|det Ad_h|^{-1}$. 
These data  define the $*$-algebra structure on $\sD(X\times H)$ by:
$$(F*G)(x,h):=\int_H \lambda^2(h') F(x,h') G(x h', h'^{-1} h)\,,\,(F^*)(x,h):=\delta(h)^{-1}\overline{F(x h, h^{-1})}$$
The algebra $C_0(X)$ is represented faithfully on $L^2(X)$ by multiplication. The induced representation  
$\Pi$ of the $*$-algebra $\sD(X\times H)$ in $L^2(X\times H)$ is given by:
$\Pi(F) (G\Psi_0)=: (\pi(F)G) \Psi_0\,,\,G\in\sD(X\times H)$, and
$$(\Pi(F)G)(x,h):=\int \lambda^2(k) F(x h^{-1},k) G(x, k^{-1} h),$$ 
where $\Psi_0=\nu_0\mt\lambda$ for some real, nonvanishing half density $\nu_0$ on $X$.  
The reduced crossed product  $C_0(X)\times_{r,\alpha} H$ is the completion of the $*$-algebra $\sD(X\times H)$ 
in the the norm coming from this representation. The canonical morphisms $i_H\in Mor(C^*(H),C_0(X)\times_\alpha H)$ 
and $i_A\in Mor(C_0(X),C_0(X)\times_\alpha H)$ are given by (the extension of):
$$(i_H(g)F)(x,h):=F(xg,g^{-1}h)\,,\,\,(i_A(f)F)(x,h):=f(x) F(x,h)\,,\,f\in \sD(X)$$

On the other hand we have  the transformation groupoid $\Delta:=X\times H$, the $*$-algebra $\sA(\Delta)$ and the $C^*$-algebra 
$C_r^*(\Delta)$. Let us choose $\om_0:=\lambda_0\mt\rho_0$, where $\lambda_0$ is defined by 
$\lambda_0(x, e)(v):=\lambda(e)(v)$ and let us define $\Psi:=\nu\mt\rho_0$. 
With such a choice the formulae for operations in $\sA(\Delta)$ and the representation $\pi_{id}$ are as follows:
$$(f*g)(x,h):=\int \lambda^2(k) f(x,k) g(x k, k^{-1} h)\,\,,\, f^*(x,h):=\overline{f(xh, h^{-1})},\,\,\pi_{id}(f)(g \Psi)=(f*g)\Psi$$
For $h\in H$ let $B_h$ denotes the operator acting on $\sA(\Delta)$ defined by a bisection $\{(x,h) : x\in X\}$.
Now define $\varphi: \sA(\Gamma)\rightarrow \sD(X\times H)$ by $(\varphi(f\om_0))(x,h):=\delta(h)^{-1/2}f(x,h)$  and let
$V: L^2(X\times H)\rightarrow L^2(X\times H)$ be a unitary defined as a push-forward by a diffeomorphism $(x,h)\mapsto (xh^{-1},h)$.
Using the definitions  above one proves:
\begin{lem} The mapping $\varphi$ is an isomorphism of  $*$-algebras. For   $h\in H\,,\, k\in\sD(X)$ and $ \om\in\sA(\Delta)$:
$i_H(h)\varphi(\om)=\varphi (B_h(\om))$ and $i_A(k)\varphi(\om)=\varphi (k \om)$ and $V \Pi(\varphi(\om))V^*=\pi_{id}(\om)$.\dowl
\end{lem}
Therefore $\varphi$ extends to an isomorphism of  $C_r^*(\Delta)$ and $C_0(X)\times_{r,\alpha} H$. 
Since $\R_*$ is abelian the reduced and universal crossed product coincide and 
we have $C^*_r(\Gamma)=C_0(\R)\times_{\alpha}\R_*$, where $(\alpha_cf)(x):=f(x/c)$.
By the universality of crossed product we have the following:
\begin{lem} Let $A$ be a  $C^*$-algebra, $G$-locally compact group and $(A,G,\alpha)$ be a dynamical system. 
Let $B:=A\times_{\alpha} G$ be the corresponding crossed product with canonical 
morphisms $i_A\in Mor(A,B)$ and $i_G\in Mor(C^*(G),B)$.
If  $A$ is generated by $X_1,X_2,\dots X_k$ and $C^*(G)$ by $Y_1,\dots,Y_n$, then $B$ is generated by 
$i_A(X_1),\dots, i_A(X_k),$ and  $i_G(Y_1),\dots, i_G(Y_n)$.\dowl
\end{lem}
Now it is clear that $C^*(\R_*)$ is generated by the the generator of one parameter group $\R\ni t\mapsto e^t\in \R_+$ acting
on $L^2(\R_*)$  and the unitary corresponding to element $-1\in\R_*$. Also $C_0(\R)$ is generated by the function 
$\R\ni x\mapsto x\in\R$. Using two previous lemmas one gets the proof of the proposition.\dow

In the remaining part of this section we express the twist $\widehat{\sT}$ as a function  of generators and compute the action 
of the comultiplication on generators. It turns out that it is given by  formulae (\ref{relkom}).

For the  morphism $\delta_0$ defined in (\ref{delta0})  we have :
$$\delta_0(B_t)=B_t\times B_t\,\,,\,\,\delta_0(J)=J\times J\,,\,\,\widehat{\delta_0}(Y)=I\mt Y+Y\mt I$$
so 
\begin{equation}\label{Delta0} \Delta_0(\hat{X})=\hat{X}\mt I+I\mt \hat{X}\,,\,\Delta_0(\hat{J})=\hat{J}\mt \hat{J}\,,
\,\,\Delta_0(\hat{Y})=I\mt \hat{Y}+\hat{Y}\mt I\end{equation}
The expressions for $\Delta_0(\hat{X})$ and $\Delta_0(\hat{Y})$ should be read in the following sense: 
since $\hat{X}$ and $\hat{Y}$ are essentialy selfadjoint on $\sD(\Gamma)$ the sums are essentialy selfadjoint on
$\sD(\Gamma)\mt\sD(\Gamma)$  and their closures are equal to  left hand sides.

For $t\in\R$ let us define sets:
\begin{equation}
\sT_t:=\{(z_1,|1+z_2|^{-t}, z_2,1) : 1+z_2\neq 0\}\,,\,\,\,,\,\,
K:=\{(z_1,sgn(1+z_2),z_2,1) : 1+z_2\neq 0\}.
\end{equation}
 Both of these sets are sections of the right and 
left projections over $E\times E'$, therefore they define (by the left multiplication) diffeomorphisms 
of $\Gamma\times U$, where $U:=\{(z,c) : z\in E'\}$,  and corresponding mappings of $\sA(\Gamma\times U)$, which  will be denoted 
by $T_t$ and $K$. Explicitely : 
$$T_t(z_1,c_1,z_2,c_2)=(z_1 |1+z_2|^t,c_1 |1+z_2|^t, z_2,c_2),$$ 
$$K(z_1,c_1,z_2,c_2)=(z_1 sgn(1+z_2),c_1 sgn(1+z_2), z_2,c_2)$$ 
and for $f\in\sA(\Gamma\times U)$: 
$$(T_t f)(z_1,c_1,z_2,c_2)=f(T_{-t}(z_1,c_1,z_2,c_2)),\,(Kf)(z_1,c_1,z_2,c_2)=f(K(z_1,c_1,z_2,c_2))$$
The formulae above define also unitary operators $\widehat{\sT}_t$ and $\widehat{K}$ on $L^2(\Gamma\times\Gamma)$.
\begin{lem} We have the equality of sets: $\sT_t\sT_r=\sT_{t+r}\,,\,\sT_1K=K \sT_1=T$. The family $\widehat{\sT}_t$ is  
a strictly continuous group of unitary multipliers of $C^*_r(\Gamma)\mt C^*_r(\Gamma)$. Let $Z$ be its generator,  then on 
$\sA(\Gamma\times U)$: $Z= \hat{X}\mt \log|\hat{Y}+I|$ and $\sD(\Gamma\times U)$ is a core for $Z$.
\end{lem}
{\em Proof: } First statement is straighforward. Also the fact that $\widehat{\sT}_t$ is a group of unitary multipliers is clear. 
So only strict continuity and a statement about $Z$ require a  proof. But this can be done exactly as in lemma \ref{op-X}. \dowl

Recall that for  a group $G$ a {\em bicharacter} $c$ is a map $c:G\times G\rightarrow S^1$ with the properties 
$c(g_1,g_2 g_3)=c(g_1,g_2) c(g_1,g_3)$ and $c(g_1 g_2, g_3)=c(g_1,g_3) c(g_2,g_3)$.
The function $Ch:\Z_2\times\Z_2\rightarrow \R$ defined by 
$$Ch(\epsilon_1,\epsilon_2):= 
\left\{\begin{array}{rcr} -1 & \,{\rm if }\, & \epsilon_1=\epsilon_2=-1\\1 &  & {\rm otherwise} \end{array}\right.$$
 is a bicharacter  and the function $\chi: (\R\times\Z_2)\times(\R\times\Z_2)\rightarrow \C$
$$\chi(x,\epsilon_1,y,\epsilon_2):=\exp(i x y)\, Ch(\epsilon_1,\epsilon_2)$$
is a bicharacter.

Let  us define the unitary $$V:=\chi(I\mt \log|\hat{Y}+I| ,I\mt sgn(\hat{Y}+I),   \hat{X}\mt I, \hat{J}\mt I),$$
 i.e.  
$V=\exp(i \hat{X}\mt \log|\hat{Y}+I|)\, Ch(I\mt sgn(\hat{Y}+I) , \hat{J}\mt I)$.
\begin{lem}
$\widehat{\sT}=V$
\end{lem}
{\em Proof:} By the previous lemma $\widehat{\sT}=\widehat{\sT}_1 \widehat{K}=\exp(\hat{X}\mt \log|\hat{Y}+I|) 
\widehat{K}$ so it remains to prove that $\hat{K}=Ch(I\mt sgn(\hat{Y}+I) , \hat{J}\mt I)$.
If  $A,B$ are  commuting operators with spectrum contained in $\{-1,1\}$ then easy  computations show that 
$Ch(A,B)=\frac12(I+A+B-AB)$.  Applying this formula to operators $A:=I\mt sgn(\hat{Y}+1)$ 
and $B:=\hat{J}\mt I$ one obtains:
$$Ch(I\mt sgn(\hat{Y}+I),\hat{J}\mt I)=\frac12(I\mt (I+sgn(I+\hat{Y}))+\hat{J}\mt (I-sgn(I+\hat{Y})))$$
and checks that this is exactly $\widehat{K}$.\dowl

In this way we obtain formulae for comultiplication on generators:
$$\Delta(\hat{Y})=V(\hat{Y}\mt I+I\mt \hat{Y})V^*,$$ $$\Delta(\hat{X})=V(\hat{X}\mt I+I\mt \hat{X})V^*,$$
$$\Delta(\hat{J})=V(\hat{J}\mt \hat{J})V^*$$
For the meaning of sums in parentheses the remark after formula (\ref{Delta0}) applies. On  $\sA(\Gamma\times U)$ we have:
$$\Delta(\hat{X})=\hat{X}\mt(\hat{Y}+I)^{-1} + I\mt\hat{X}\,,\,\,\Delta(\hat{Y})=I\mt\hat{Y}+\hat{Y}\mt I+\hat{Y}\mt\hat{Y}$$
and these are relations (\ref{relkom}).

%%%%%%%%%%%%%%%%%%%%%%%%%%%%%%%%%%%%%%%%%%%%%%%%%%%%%%%%%%%%%%%%%%%%%%%%%%%
%%%%%%%%%%%%%%%%%%%%%%%%%%%%%%%%%%%%%%%%%%%%%%%%%%%%%%%%%%%%%%%%%%%%%%%%%%%
%%%%%%%%%%%%%%%%%%%%%%%%%%%%%%%%%%%%%%%%%%%%%%%%%%%%%%%%%%%%%%%%%%%%%%%%%%%%
\section{Poisson-Lie strucure}
\newcommand{\bGamma}{\mathbf{G}}
In this section we consider the family of groupoids $\Gamma_s\,,s\in\R$ over $A$ defined in the 
section (\ref{setup}). We define 
$\lambda_s(0,a)(\partial_b+s\partial_a):=1$. 
The  corresponding left and right invariant half densities
on $\Gamma_s$ are given by:
$$\lambda_s(b,a)(\partial_b+s \partial_a)=\rho_s(b,a)(\partial_b+s \frac{a}{1+s b}\partial_a):=
|1+s b|^{-1/2}$$
We put $\om_s:=\lambda_s\mt\rho_s$ and identify $\sA(\Gamma_s)$ with $\sD(\Gamma_s)$ with 
multiplication and involution defined  by $(f\om_s)(g\om_s)=:(f *_s g)\om_s$  and $(f\om_s)^*=:(f^{*_s})\om_s$:
\notka{prod-s}\begin{align}\label{prod-s}
(f *_s g)(b,a):=\int \frac{{\rm d  c}}{|1+s c|} f(c, a+s (c-b))\,g( \frac{b-c}{1+s c},\frac{a}{1+s c})=\nonumber\\
=\int \frac{{\rm d  c}}{|1+s c|} f(\frac{b-c}{1+s c}, a-s c \frac{1+s b}{1+s c})\,g( c,a \frac{1+s c}{1+s b})
\end{align}
and $f^{*_s}(b,a):=\overline{f(\frac{-b}{1+s b}, \frac{a}{1+s b})}$. We will write $f*g$ and $f^*$ 
instead of $f*_0 g$ and $f^{*_0}$.  
For $M>1$ let $K_M:=\{(b,a)\in\Gamma_0: |b|\leq M\,,\,\frac1M\leq |a|\leq M\}$, it is clear that 
any $f\in\sD(\Gamma_0)$ is supported in $K_M$ for sufficiently large $M$. 
The product $*_s$ is in fact defined for all $f,g\in\sD(\Gamma_0)$:
\begin{lem}\label{lem-prod-s}\notka{lem-prod-s} Let $f,g\in\sD(\Gamma_0)$ have supports in $K_M$. 
Then for any $s\in\R$ the function $f*_sg$ 
defined by (\ref{prod-s}) is  smooth and has the support in the a set 
$\{(b,a)\in\Gamma_0: |b| \leq M(2+|s|M)\,,\,\frac{1}{M(1+M^2 |s|)}\leq |a|\leq M(1+|s|M)\}$. 
In particular, for $|s|\leq \delta$ functions $f*_s g$ have supports contained in a fixed compact set.
\end{lem}
{\em Proof:} Let $(b,a,s)\in\Gamma_0\times\R$. It is straightforward that the function 
$$\R\ni b'\mapsto f(b',a+s(b'-b)) g(\frac{b-b'}{1+s b'},\frac{a}{1+s b'})$$ 
is smooth and has compact support contained  in the set $\{b'\in\R: 1+s b'\neq 0\}$, so 
smoothness follows. For $(f*_s g)(b,a)\neq 0$ it is necessary that there exists
$b'$ such that 
$$|b'|\leq M\,,\,\frac1M\leq |a+s(b'-b)|\leq M\,,\,\left|\frac{b-b'}{1+s b'}\right|\leq M\,,\,
\frac1M\leq \left|\frac{a}{1+s b'}\right|\leq M.$$
If $|a|< \frac1{M(1+M^2 |s|)}$  then
$$\frac1{M}\leq |a+s(b'-b)|\leq |a|+|s||b'-b| <  \frac{1}{M(1+M^2 |s|)}+|s| M |1+s b'|\leq$$
$$\leq \frac{1}{M(1+M^2 |s|)}+ |s|M^2 |a|<\frac{1}{M(1+M^2 |s|)}+|s|M^2 \frac{1}{M(1+M^2 |s|)}=\frac1M$$
In a similar way, if $|a|>M(1+|s|M)$ then:
$$1+|s|M\geq 1+|s||b'|\geq |1+s b'|\geq \frac{|a|}{M}>1+|s|M$$
The estimates for $|a|$ are proven. If $|b|>M(2+|s|M)$ then 
$$M(2+|s|M)<|b|\leq |b-b'|+|b'|\leq M|1+s b'|+M\leq M(1+|s|M)+ M=M(2+ |s|M).$$
\dowl

The choice  $\nu(a)(\partial_a):=\frac{1}{|a|}$  defines real, 
non vanishing half density  $\Psi_s:=\rho_s\mt\nu$ on  $\Gamma_s$ 
and short calculation gives: $\Psi_s(b,a)(\partial_b,\partial_a)=:\Psi_0(b,a)(\partial_b,\partial_a)= \frac{1}{|a|}$. 
This makes possible identification of  all spaces 
$L^2(\Gamma_s)$ with $L^2(\Gamma_0)$. The identity representation of $\sD(\Gamma_s)$ is then given by:
$\pi_s(f)(g\Psi_0)=(f *_s g) \Psi_0$ for $f,g\in\sD(\Gamma_s)$.
The norms on $\sD(\Gamma_s)$ defined by $\om_s$ will be denoted by $||f||_{l,s}$, 
$||f||_{r,s}$, $||f||_{0,s}$ and $||f||_s$ 
is the norm of $\pi_s(f)$. Finally let us define for $f,g\in\sD(\Gamma_0)$:
$$\{f,g\}:=(a-1) [(\partial_a f)*(b g)-(\partial_a g)*(b f)].$$
In this formula $(a-1)f$ and $bf$ denote functions $((a-1)f)(b,a):=(a-1)f(b,a)$ and $(bf)(b,a):=bf(b,a)$.
\begin{lem}\label{poisson}
$\{\cdot,\cdot\}$ is a Poisson bracket on $\sD(\Gamma_0)$ and $\{f_1^*, f_2^*\}=\{f_2,f_1\}^*$.
\end{lem}
{\em Proof:} The mappings $f\mapsto (a-1)\partial_a f$ and $f\mapsto bf$ are commuting derivations of a commutative algebra 
$(\sD(\Gamma_0),*)$, moreover $(a-1)\partial_a(f^*)=((a-1)\partial_a f)^*$ and $bf^*=-(bf)^*$\dowl

It is clear that if  $f\in\sD(\Gamma_0)$, then there exists $\delta>0$ such that 
$f\in\sD(\Gamma_s)$  for all $|s|\leq\delta$, e.g.
if support  of $f$ is contained in a set $K_M$,  take  any $\delta<\frac1{M^2}$. 
Let us define linear spaces $$D(Q_s):=\{f\in\sD(\Gamma_0) : |\delta|\leq |s| \Rightarrow f\in\sD(\Gamma_\delta)\}.$$
This family of subspaces has properties:
$$|r|\leq|s|\Rightarrow D(Q_s)\subset D(Q_r)\,, \,D(Q_s)^*=D(Q_s)\,,\,\bigcup_{s\neq0}D(Q_s)=\sD(\Gamma_0)$$
Let us define linear map $Q_s: \sD(\Gamma_0)\supset D(Q_s)\ni f\mapsto f\in\sD(\Gamma_s)$. 
With this definition we have
\begin{prop} \label{quant}
\begin{enumerate}
\item $Q_0=id$;
\item $\lim_{s\rightarrow 0} ||Q_s(f^*)-Q_s(f)^{*_s}||_s=0$;
\item $\lim_{s\rightarrow 0} ||Q_s(f)*_s Q_s(g)-Q_s(f*g)||_s=0$;
\item $\lim_{s\rightarrow 0} ||\frac1s[Q_s(f),Q_s(g)]-Q_s(\{f,g\})||_s=0$;
\item The function $s\mapsto ||Q_s(f)||_s$ is continuous for $s\neq 0$ and 
lower semicontinuous at $s=0$.
\end{enumerate}
\end{prop}
{\em Proof:} First statement is evident. To prove statements (2)-(4) it is enough to prove that 
norms $||\cdot||_{0,s}$ converge to $0$. In fact we claim that convergence of norms $||\cdot||_{l,s}$ is sufficient.
In the following to simplify notation we write $Q_s(f) Q_s(g)$ instead of $Q_s(f)*_sQ_s(g)$ 
and $Q_s(f)^*$ instead of $Q_s(f)^{*_s}$. Let us compute:
$$||Q_s(f^*)^*-Q_s(f)||_{r,s}=||Q_s(f^*)-Q_s(f)^*||_{l,s}=
||Q_s((f^*)^*)^*-Q_s(f^*)||_{l,s}=||Q_s(f_1^*)^*-Q_s(f_1)||_{l,s},$$
where $f_1:=f^*$. So once we know that $\lim_{s\rightarrow 0} ||Q_s(f^*)^*-Q_s(f)||_{l,s}=0$ 
for any $f\in\sD(\Gamma_0)$, the second statement is proven.\\
In a similar way:
$$||Q_s(f) Q_s(g)-Q_s(f*g)||_{r,s}=||Q_s(g)^* Q_s(f)^*-Q_s(f*g)^*||_{l,s}$$

Let  $f_1:=f^*\,,\,g_1:=g^*$, then 
$$Q_s(g)^* Q_s(f)^*-Q_s(f*g)^*=Q_s(g_1^*)^*Q_s(f_1^*)^*-Q_s(f_1^**g_1^*)^*=$$
$$\left(Q_s(g_1^*)^*-Q_s(g_1)\right) \left(Q_s(f_1^*)^*-Q_s(f_1)\right)+ \left(Q_s(g_1^*)^*-Q_s(g_1)\right) Q_s(f_1)   + 
Q_s(g_1)\left(Q_s(f_1^*)^*-Q_s(f_1)\right)+$$
$$+ \left(Q_s(g_1*f_1)-Q_s((g_1*f_1)^*)^*\right)+Q_s(g_1) Q_s(f_1)-Q_s(g_1 *f_1)$$
And once more:
$$||[Q_s(f),Q_s(g)]-s Q_s(\{f,g\})||_{r,s}=||[Q_s(f),Q_s(g)]^*-s Q_s(\{f,g\})^*||_{l,s}$$
$$[Q_s(f),Q_s(g)]^*-s Q_s(\{f,g\})^*=Q_s(g)^*Q_s(f)^*-Q_s(f)^*Q_s(g)^*-s Q_s(\{f,g\})^*=$$
$$=Q_s(g_1^*)^*Q_s(f_1^*)^*-Q_s(f_1^*)^*Q_s(g_1^*)^*-s Q_s(\{f_1^*,g_1^*\})^*=[Q_s(g_1^*)^*,Q_s(f_1^*)^*]-sQ_s(\{g_1,f_1\}^*)^*=$$
$$=[Q_s(g_1^*)^*-Q_s(g_1),Q_s(f_1^*)^*-Q_s(f_1)]+(Q_s(g_1^*)^*-Q_s(g_1)) Q_s(f_1)+[Q_s(g_1),Q_s(f_1^*)^*-Q_s(f_1)]+$$
$$-s (Q_s(\{g_1,f_1\}^*)^*-Q_s(\{g_1,f_1\}))+ ([Q_s(g_1),Q_s(f_1)]-s Q_s(\{g_1,f_1\}))$$

Using the equalities above one can see that to prove statements 2-4 of the proposition it is enough to prove 
\begin{lem}\label{granica-s} For any $f\in\sD(\Gamma_0)$ there is $\delta>0\,,\,C>0$ such that $||Q_s(f)||_{l,s}\leq C$ 
for $|s|\leq \delta$; Moreover for any $f,g\in\sD(\Gamma_0)$: 
$$\lim_{s\rightarrow 0} ||Q_s(f^*)^{*_s}-Q_s(f)||_{l,s}= 
\lim_{s\rightarrow 0} ||Q_s(f)*_s Q_s(g)-Q_s(f*g)||_{l,s}=$$
$$=\lim_{s\rightarrow 0} ||\frac1s[Q_s(f),Q_s(g)]-Q_s(\{f,g\})||_{l,s}=0.$$
\end{lem}
{\em Proof:}  Let $supp f\subset K_M\subset \Gamma_s$ for $|s|\leq \delta$. Then 
$||Q_s(f)||_{l,s}=\sup_{a\in A}\int_{B_s} \frac{db'}{|1+s b'|}|f(b',a+s b')|$. For $\delta<\frac1M$: 
$\frac{1}{|1+s b'|}\leq \frac1{1-|s| M}$ on the support of $f$, 
therefore $||Q_s(f)||_{l,s}\leq \frac{2M}{1-|s| M} \sup |f|$ 
and the estimate proves the first statement.

Let $supp f\cup supp f^* \subset K_M$, for any   $|s| < \frac1{M^2}$   $supp f,\,supp f^* \subset\Gamma_s$. 
Define   $F(b,a):=(Q_s(f^*)^{*_s}-Q_s(f))(b,a)=f(\frac{b}{1+s b},\frac{a}{1+s b})-f(b,a)$.
%$||F||_{l,s}=\sup_{a\in A}\int_{B_s} \frac{db'}{|1+s b'|}|F(b',a+s b')|$. 
Fix $a\in A$ and consider the function 
$b'\mapsto |F(b',a+s b')|$. Its support  is contained in $S_1\cup S_2$ for 
$S_1:=\{b': \left|\frac{b'}{1+s b'}\right|\leq M\,,\,\frac1M\leq \left|\frac{a+s b'}{1+s b'}\right|\leq M\}$ and 
$S_2:=\{b': |b'|\leq M\,,\,\frac1M\leq|a+s b'|\leq M\}$. Now, for $|s| \frac1{2M}$ we 
get that union of these sets
is contained in the set $S_3:=\{b': |b'|\leq 2 M\}$ and on this set we have estimates: 
$$\left|\frac1{1+s b'}\right|\leq 2\,\,{\rm and} \,\,\frac1M\leq \left|\frac{a+s b'}{1+s b'}\right|\leq 2 M.$$
Using these estimates and equalities: $\frac{b'}{1+s b'}=b'+\frac{- s b'^2}{1+s b'}$, 
$\frac{a+s b'}{1+s b'}=a+s b'+\frac{-s b'(a+s b')}{1+s b'}$ we get:
$$|F(b',a+s b')|\leq \sup ||f'||\left(\left|\frac{- s b'^2}{1+s b'}\right|+\left|\frac{-s b'(a+s b')}{1+s b'}\right|\right)\leq 
|s| \sup ||f'||( 2 |b'|^2 + 2 M|b'|)$$
It follows that $\lim_{s\rightarrow 0} ||F||_{l,s}=0$, so $\lim_{s\rightarrow 0} ||Q_s(f^*)^{*_s}-Q_s(f)||_{l,s}=0.$

For sufficiently small $|s|$ the functions $(Q_s(f)Q_s(g)-Q_s(f*g)$ have support in a fixed compact subset of $\Gamma_s$, 
so we have to show that $Q_s(f)Q_s(g)-Q_s(f*g)$ converges to $0$ uniformly. So choose $M$ such that 
$supp\,f\cup supp\,g \cup supp (f*g)\subset K_M\,, |s|M^2<1$. Then 
$$(Q_s(f)Q_s(g)-Q_s(f*g))(b,a)=\int_{-M}^{M} {\rm d b'}\left(\frac{1}{1+s b'}f(b',a+s(b-b')) 
g(\frac{b-b'}{1+s b'},\frac{a}{1+s b'})- f(b',a) g(b-b',a)\right).$$ 
For $|b'|\leq M$ we have 
$$\left|\frac{1}{1+s b'}f(b',a+s(b-b')) g(\frac{b-b'}{1+s b'},\frac{a}{1+s b'})-f(b',a) g(b-b',a)\right|\leq$$
$$\leq (1+s b')^{-1}\left|f(b',a+s(b-b')) g(\frac{b-b'}{1+s b'},\frac{a}{1+s b'})-f(b',a) g(b-b',a)-s b'f(b',a) g(b-b',a)\right|
\leq$$
$$\leq (1-|s| M)^{-1}\left(\left|f(b',a+s(b-b')) g(\frac{b-b'}{1+s b'},\frac{a}{1+s b'})-f(b',a) g(b-b',a)\right|+
|s|M\sup|f g|\right)$$
For the first term there is an estimate:
$$\left|f(b',a+s(b-b')) g(\frac{b-b'}{1+s b'},\frac{a}{1+s b'})-f(b',a) g(b-b',a)\right|\leq$$
$$\leq |s| \sup|g|\sup|\partial_af||b-b'|+\sup|f|\left|g(\frac{b-b'}{1+s b'},\frac{a}{1+s b'})-g(b-b',a)\right|$$
and finally: 
$$\left|g(\frac{b-b'}{1+s b'},\frac{a}{1+s b'})-g(b-b',a)\right|\leq 
|s|M\sup||g'||\left(\left|\frac{b-b'}{1+sb'}\right|+
\left|\frac{a}{1+sb'}\right|\right)\leq$$
$$\leq |s|M\sup||g'||(1-|s|M)^{-1}(|b|+|a|+M),$$
since $b,a$ are in a fixed compact set,  convergence is uniform 
and $\lim_{s\rightarrow 0} ||Q_s(f)*_s Q_s(g)-Q_s(f*g)||_{l,s}=0$.

In a similar way, to prove the third equality we have to show that the functions 
$\frac1s[Q_s(f),Q_s(g)]-Q_s(\{f,g\})$ converge to $0$ uniformly. Again, let $f,g,\{f,g\}$ 
be supported in $K_M \subset\Gamma_s$:
$$(Q_s(f)Q_s(g))(b,a)=\int_{-M}^M \frac{{\rm d b'}}{1+s b'} f(b', a+s (b'-b))\,
g( \frac{b-b'}{1+s b'},\frac{a}{1+s b'})$$
$$(Q_s(g)Q_s(f))(b,a)=\int_{-M}^M \frac{{\rm d b'}}{1+s b'} g(\frac{b-b'}{1+s b'}, a-sb'\frac{1+s b}{1+s b'})
\,f( b', a\frac{1+sb'}{1+s b})$$
$$(Q_s(\{f,g\}))(b,a)=(a-1)\int_{-M}^{M} {\rm d b'}
\left[(\partial_af)(b',a)\, (b-b')\,g(b-b',a)-b'\,f(b',a)\, (\partial_ag)(b-b',a)\right]$$
(in the formula for $Q_s(g)Q_s(f)$ we use second expression of (\ref{prod-s})). 
This computation is straightforward and can be done as in the previous point. This ends the prove of the lemma and statements 
(2)-(4) of the proposition (\ref{quant}) \dowl
%\\\qed

Now we come to the last point. Let's start with lower semicontuity at $s=0$. So we have to prove that:
$\forall\,f\in\sD(\Gamma_0)\, \forall\,\epsilon>0\,\exists\,\delta>0\,
\left(|s|<\delta\, \Rightarrow ||Q_s(f)||_s> ||f||-\epsilon\right)$ 
( $||f||$--denotes the operator norm of $\pi_0(f)$ on $L^2(\Gamma_0)$). The following simple lemma reduces the problem to 
strong convergence of operators $Q_s(f)$:
\begin{lem} Let $H$ be a Hilbert space, $A\in B(H)$,  
$A_s$ - family of bounded operators defined for some neighbourhood of $0\in\R$ and $V\subset H$ a dense subspace.
Assume that for every $v\in V$ $\lim_{s\rightarrow 0} A_s v=Av$. Then $ \forall\,\epsilon>0\,\exists\,\delta>0\,
\left(|s|<\delta\, \Rightarrow ||A_s||> ||A||-\epsilon\right)$
\dowl\end{lem}
So we have to prove that for $f,g\in\sD(\Gamma_0)$ $f*_s g -f*g$ converges to $0$ in $L^2(\Gamma_0)$ as $s$ goes to $0$. 
But for sufficiently small $|s|$ these functions have supports in a fixed compact set, so it is enough to prove uniform 
convergence and this can be done as in the previous points.

Now, for $s,r\neq 0$, consider the map  $\Phi_{sr}:\Gamma_r\ni(b,a)\mapsto  (\frac{r}{s}b,a)\in\Gamma_s$. This is a 
diffeomorphism and an isomorphism of groupoids. It defines a $*$-isomorphism $\Phi_{sr}: \sD(\Gamma_r)\lra \sD(\Gamma_s)$ and
unitary operator $V_{sr}$ on $L^2(\Gamma_0)$. The formulae are:
$$(\Phi_{sr}f)(b,a):=|\frac{s}{r}|f(\frac{s}{r} b,a)\,\,,\,
V_{sr}(f\Psi_0)=\sqrt{|\frac{r}{s}|}(\Phi_{sr}f)\Psi_0\,,\,f\in\sD(\Gamma_r)$$
Moreover, we have: $$V_{rs}V_{sr}=I\,,\,V_{rs}\pi_s(f)V_{sr}=\pi_r(\Phi_{rs}(f))\,,\,f\in\sD(\Gamma_s)$$

Let $f\in\sD(\Gamma_s)$ then there exists $|s|>\delta>0$ such that $f\in\sD(\Gamma_{s+\epsilon})$ for $|\epsilon|<\delta$.
$$||Q_{s+\epsilon}(f)||_{s+\epsilon}=||\pi_{s+\epsilon}(f)||=||V_{s s+\epsilon}\pi_{s+\epsilon}(f)V_{ s+\epsilon s}||=
||\pi_s(\Phi_{s s+\epsilon}(f))||=||Q_s(\Phi_{s s+\epsilon}(f))||_s$$
Now, it is straightforward to check that the functions 
$\Phi_{s s+\epsilon}(f)-f$ have support in a fixed compact set and tend to $0$
uniformly as $\epsilon$ goes to $0$. This  ends the proof of the  proposition (\ref{quant}).
\dow

The groupoid $\Gamma_0$ is a (trivial) bundle of groups $\R$, as such it can be identified with its Lie algebroid 
$\R\times A$. 
For $f\in\sD(\Gamma_0)$ let  $(\sF f)(\beta,a):=\int db \,e^{-2\pi i \beta b}\, f(b,a)$ be the (partial) Fourier transform. 
Let $\sB:=\sF(\sD(\Gamma_0))$, then it is a $*$-subalgebra of smooth functions (with  respect to pointwise multiplication 
and complex conjugation as an involution) on a dual bundle. For $F,G\in\sB$ let us define:
$$\{F,G\}:=-i\sF(\{\sF^{-1}F,\sF^{-1}G\})$$
Then straightforward calculation gives:
\begin{equation}\label{nawias} \{F,G\}(\beta, a)=\frac{a-1}{2 \pi}\left((\partial_aF)(\partial_\beta G)(\beta,a)-
(\partial_a G)(\partial_\beta F)(\beta,a)\right)\end{equation}

On the other hand recall the morphism $m_B^T: \Gamma_0\rel \Gamma_0\times\Gamma_0$ defined in section 2
$$Gr(m_B^T):=\{(b_2a_1, a_1, b_2, a_2; b_2 a_1, a_1 a_2)\}$$
Applying the cotangent lift (see \cite{SZ2}) we obtain the morphism (of symplectic groupoids) $T^*(\Gamma_0)$ and 
$T^*(\Gamma_0)\times T^*(\Gamma_0)$, its base map is nothing but the multiplication in the subgroup 
$(TA)^0\subset T^*(\Gamma_0)$. The map $\R\times A\ni(\beta,a)\mapsto \beta db(0,a)\in (TA)^0$ identifies the group $(TA)^0$ with
$\R\times A$ with the multiplication:
\begin{equation}\label{ax-n}
(\beta_1,a_1)(\beta_2,a_2):=(\beta_1+a_1^{-1}\beta_2, a_1 a_2)
\end{equation} (which is again 'ax+b' 
group in different presentation). One easily checks that the bracket (\ref{nawias}) is a Poisson-Lie bracket on this group.

Finally let us consider a comultiplication. The groupoid $\Gamma_0$ is related to the double Lie group $(G;A,B)$, so $C^*_r(\Gamma_0)$ 
is equipped with the  comultiplication $\Delta_0$, which is an extension of the mapping 
$\widehat{\delta_0}$ given, for $f\in\sD(\Gamma_0)$ and  $F\in\sD(\Gamma_0\times\Gamma_0)$, by (compare (\ref{delta0-form})):
\begin{equation}\label{ddelta0}
(\widehat{\delta_0}(f)F)(b_1,a_1,b_2,a_2)=\int d b\, f\left(b,a_1 a_2\right)\,F(b_1-b, a_1, -\frac{b}{a_1}+b_2, a_2)
\end{equation}
For $\tilde{f}=\sF{f}\in\sB$ and $\tilde{F}=(\sF\times \sF)(F)\in\sB\mt\sB$ let us define:
$$\Delta(\tilde{f})\tilde{F}:=(\sF\times \sF)(\widehat{\delta_0}(f)F)$$
Routine calculations show that:
$$[\Delta(\tilde{f})\tilde{F}](\beta_1,a_1,\beta_2,a_2)=\tilde{f}(\beta_1+\frac{\beta_2}{a_1},a_1 a_2)\tilde{F}(\beta_1,a_1, \beta_2,a_2)$$
so $\Delta$ is a comultiplication given by the group structure (\ref{ax-n}).

Each of groupoids  $\Gamma_s$ comes together with a relation $\delta_s:=\tilde{m}_C^T$ (the relation $\tilde{m}_C$ is defined in (\ref{mc})).
Although $\delta_s$ (for $s\neq 0$) is {\em not} a morphism $\Gamma_s\rel \Gamma_s\times \Gamma_s$, one can define the mapping 
$\widehat{\delta_s}$ (in fact this is the restriction of the mapping $\hat{\delta}$ from the section \ref{stwist}):
\begin{eqnarray}\label{ddelta-s}(\widehat{\delta_s}(f)F)(b_1,a_1,b_2,a_2)=
\int \frac{d b}{|1+s b|} f\left(b,s b+(a_1-s b_1) (a_2-s b_2)\right)\times\\
\nonumber 
\times F\left(\frac{b_1-b}{1+ s b}, \frac{a_1}{1+s b}, \frac{-b+b_2 (a_1-s b_1)}{a_1-s b_1+s b}, \frac{a_2(a_1- s b_1)}{a_1-s b_1+s b}\right)
\end{eqnarray}
For $s=0$ this is (\ref{ddelta0}). In general $\widehat{\delta_s}(f)F$ (for $s\neq 0$) is not an element 
of $\sD(\Gamma_0\times\Gamma_0)$, this is the case for sufficiently small $|s|$, hewever:
\begin{lem}\label{supp-delta-s}
Let $f,F$ have supports in $K_M$ and $K_M\times K_M$ respectively ($K_M$ was defined right after formula (\ref{prod-s})).
There exists $k>1$ such that for $|s|<\frac{1}{k M^2}$:
$$(b_1,a_1,b_2,a_2)\in supp(\widehat{\delta_s}(f)F) \Rightarrow $$ 
$$\frac1M-\frac{1}{kM^2}\leq |a_1|\leq M+\frac1k\,,\,|b_1|\leq 2M+\frac1k\,,\,|b_2|\leq2 M^2+M(1+\frac2k)\,,\,
\frac1M-\frac2{k M}\leq |a_2|\leq M(1+\frac2k)$$
In particular for any $M>1$ $supp(\widehat{\delta_s}(f)F)$ is contained in a fixed compact set for sufficiently small $|s|$.
\end{lem}
{\em Proof:} The proof is similar to that of lemma \ref{lem-prod-s}. One writes inequalities resulting from 
assumptions about supports of $f$ and $F$ and, after some manipulations, gets estimates as in the lemma.\dowl

\noindent
Let us end with the following:
\begin{prop}
Let $F\in\sD(\Gamma_0\times\Gamma_0)$ and $f\in\sD(\Gamma_0)$. There exists $\epsilon>0$ such that for $|s|\leq\epsilon$ 
$\widehat{\delta_s}(f)F\in \sD(\Gamma_0\times\Gamma_0)$ is supported  in a fixed compact set  and 
$$\lim_{s\rightarrow 0}(\widehat{\delta_s}(f)F-\widehat{\delta_0}(f)F)=0\,\,\,{\rm in\,} C^*_r(\Gamma_0)$$
\end{prop}
{\em Proof:} First statement follows from lemma \ref{supp-delta-s} and standard arguments for smoothness. To prove the second statement, 
one proceeds as in the proof of lemma (\ref{granica-s}): it is sufficient to prove uniform convergence and this can be easily 
done using results of lemma \ref{supp-delta-s}.\dow

%%%%%%%%%%%%%%%%%%%%%%%%%%%%%%%%%%%%%%%%%%%%%%%%%%%%%%%%%%%%%%%%%%%%%%%%%%%%%%%
%%%%%%%%%%%%%%%%%%%%%%%%%%%%%%%%%%%%%%%%%%%%%%%%%%%%%%%%%%%%%%%%%%%%%%%%%%%%%%%
\section{Appendix}

Here we collect some formulae proven in \cite{DLG} and used in this paper.
%
%$G$ is a Lie group. $A,B,C$ are closed subgroups; $\gotg, \gota,\gotb, \gotc$ are corresponding Lie algebras.
$(G;B,C)$ is a double Lie group,   $\gotg, \gotb,\gotc$  are corresponding Lie algebras and $\gotg=\gotb\oplus \gotc$ 
(direct sum of vector spaces).

\noindent{\em Modular functions.}
Let  $P_B, P_C$ %, P_C,\tilde{P}_B$ 
be projections in $\gotg$ corresponding to the decomposition $\gotg=\gotb\oplus \gotc$. %=\gotc\oplus\gotb.$$
Let us define:
\notka{modular-functions}\begin{equation}\label{modular-functions}
j_B(g):=|\det(P_B Ad(g)|_{\gotb})|\,,\,\,j_C(g):=|\det(P_C Ad(g)|_{\gotc})|
\end{equation}

\noindent {\em The choice of $\om_0$.}
Choose a real half-density $\mu_0\neq 0$ on $T_eC$ and define left-invariant half-density on $G_B$ by
$$\lambda_0(g)(v):=\mu_0(g^{-1} v)\,,\,v\in \lma T^l_g G_B$$
Then the corresponding right-invariant  half-density is given by:
%$$\rho_0(g)(w):=|\psi_B(a_L(g))|^{-1/2}\mu_0(w g^{-1})\,,\,w\in\lma T^r_g G_A$$
$$\rho_0(g)(w):=j_C(b_L(g))^{-1/2}\mu_0(w g^{-1})\,,\,w\in\lma T^r_g G_B.$$
{\em Multiplication and comultiplication in $\sA(G_B)$}
After the choice of $\om_0$ as above,  the multiplication in $\sA(G_B)$ reads:
$(f_1\om_0)(f_2\om_0)=:(f_1*f_2)\om_0$ and
\notka{mult-ga}\begin{equation}\label{mult-ga}
(f_1*f_2)(g)=\int_C d_lc\, f_1(b_L(g) c)f_2(c_L(b_L(g) c)^{-1}g)=$$
$$\int_C d_rc\, j_C(b_L(c b_R(g)))^{-1} f_1(g  c_R(c b_R(g))^{-1}) f_2(c b_R(g)),
\end{equation}
where $d_l c$ and $d_r c$ are left and right Haar measures on $C$ defined by $\mu_0$.

The $||\cdot||_l$ defined by this $\om_0$ is given by:
\notka{l-norm}\begin{equation}\label{l-norm}
||f||_l=\sup_{b\in B}\int_C d_l c\, |f(b c)|
\end{equation}
Let $\delta_0:=m_C^T: G_B\rel G_B\times G_B$. The formula for $\hat{\delta}_0$ reads
$$\hat{\delta}_0(f\om_0)(F(\om_0\otimes\om_0))=:(\hat{\delta}_0(f) F)(\om_0\otimes \om_0)$$
\notka{delta0-form}\begin{equation}\label{delta0-form}
(\hat{\delta}_0(f) F)(b_1 c_1, b_2 c_2)=
\int_C d_lc\, j_C(c_L(b_2 c))^{-1/2}f(b_1 b_2 c)F(c_L(b_1 b_2 c)^{-1} b_1 c_1, b_R(b_2 c) c^{-1} c_2)
\end{equation}
The mappings $\widehat{(id \times \delta_0)}$ and $\widehat{(\delta_0\times id)}$ are given by:
\notka{id-delta0-form}\begin{equation}\label{id-delta0-form}
[\widehat{(id \times \delta_0)}(F_1)F_2](b_1 c_1, b_2 c_2, b_3 c_3)=
\int_{C\times C} d_lc' d_lc'' F_1(b_1 c', b_2 b_3 c'')\times$$
$$\times F_2(c_L(b_1 c')^{-1} b_1 c_1, c_L(b_2 b_3 c'')^{-1} b_2 c_2, b_R(b_3
c'')^{-1} c''^{-1} c_3) j_C(c_L(b_3 c''))^{-\frac12}
\end{equation}
%%%%%%%%%
\notka{delta0-id-form}\begin{equation}\label{delta0-id-form}
[\widehat{(\delta_0\times id)}(F_1)F_2](b_1 c_1, b_2 c_2, b_3 c_3)=
\int_{C\times C} d_lc' d_lc'' F_1(b_1 b_2 c', b_3 c'')\times$$
$$\times F_2(c_L(b_1 b_2 c')^{-1} b_1 c_1, b_R(b_2 c')c'^{-1} c_2, c_L(b_3
c'')^{-1} b_3 c_3) j_C(c_L(b_2 c'))^{-\frac12}
\end{equation}

%%%%%%%%%%%%%%%%%%%%%%%%%%%%%%%%%%%%%%%%%%%%%%%%%%%%%%%%%%%%%%%%%%%%%%%%%%%%%%


\begin{thebibliography}{}
\bibitem{BS} G. Skandalis, {\em Duality for locally compact 'quantum groups'} (joint work with S. Baaj), 
Mathematisches Forschungsinstitut Oberwolfach, Taungsbericht 46/1991, $C^*$-algebren, 20,10-26.10.1991,p. 20;
\bibitem{VV} S. Vaes, L. Vainerman, {\em Extensions of locally compact quantum groups and the bicrossed product 
constraction}, Adv. in Math {\bf 175} (1) (2003), 1-101;
\bibitem{BV} S. Baaj, G. Skandalis, S. Vaes {\em Non-semi-regular quantum groups coming from number theory}, 
Comm. Math. Phys. {\bf 235} (1) (2003), 139-167;
\bibitem{SZ1} S. Zakrzewski {\em Quantum and Classical pseudogroups I}  Comm. Math. Phys. {\bf 134} (1990), 347-370;
\bibitem{SZ2} S. Zakrzewski {\em Quantum and classical pseudogroups. II. Differential and
symplectic pseudogroups},  Comm. Math. Phys., {\bf 134} (2) (1990), 371-395;
\bibitem{DLG} P. Stachura, {\em From double Lie groups to quantum groups}, Fund. Math. {\bf 188} (2005), 195-240.
\bibitem{DG}  P. Stachura, {\em Differential groupoids and $C^*$-algebras}, math.QA/9905097, 
for a short exposition see: {\em $C^*$-algebra of a differential groupoid}, Banach Center Publ 51, 
Inst. Math. Polish Acad. Sci., 2000, 263-281.
\bibitem{Wor1} S. L. Woronowicz, K. Napi\'{o}rkowski, {\em Operator theory in the C*-algebra framework}, 
 Reports on Math. Phys. {\bf 31} (3) (1992), 353-371.
\bibitem{tuset} S. Neshveyev, L. Tuset, {\em Deformation of $C^*$-algebras by cocycles on locally compact quantum groups}, arXiv:1301.4897.
\end{thebibliography}
\end{document}